\documentclass[12pt]{amsart}

\usepackage{amsmath}
\usepackage{amsfonts}
\usepackage{amssymb}
\usepackage{latexsym}
\usepackage{euscript}

\newtheorem{theorem}{Theorem}[section]
\newtheorem{proposition}[theorem]{Proposition}
\newtheorem{corollary}[theorem]{Corollary}
\newtheorem{lemma}[theorem]{Lemma}

\theoremstyle{remark}
\newtheorem{remark}[theorem]{Remark}

\theoremstyle{definition}
\newtheorem{example}[theorem]{Example}
\numberwithin{equation}{section}

\def\sD{{\mathfrak D}}      
   \def\sH{{\mathfrak H}}   
      
\def\sM{{\mathfrak M}}   \def\sN{{\mathfrak N}}   
      
   \def\sT{{\mathfrak T}}

      \def\dC{{\mathbb C}}
\def\dD{{\mathbb D}}

   \def\dN{{\mathbb N}}   
      \def\dR{{\mathbb R}}
   \def\dT{{\mathbb T}}   
      
   \def\dZ{{\mathbb Z}}

\def\cA{{\mathcal A}}   \def\cB{{\mathcal B}}   \def\cC{{\mathcal C}}
      \def\cF{{\mathcal F}}
\def\cG{{\mathcal G}}   \def\cH{{\mathcal H}}   
   \def\cK{{\mathcal K}}   \def\cL{{\mathcal L}}
\def\cM{{\mathcal M}}      
\def\cP{{\mathcal P}}   \def\cQ{{\mathcal Q}}   
\def\cS{{\mathcal S}}   \def\cT{{\mathcal T}}   \def\cU{{\mathcal U}}
\def\cV{{\mathcal V}}   \def\cW{{\mathcal W}}   \def\cX{{\mathcal X}}

\def\RE{{\rm Re\,}}
\def\IM{{\rm Im\,}}

\def\wt{\widetilde}
\def\wh{\widehat}

\def\f{\varphi}

\def\uphar{{\upharpoonright\,}}
\def\ovl{\overline}

\def\span{{\rm span\,}}
\def\ran{{\rm ran\,}}
\def\dom{{\rm dom\,}}
\def\cspan{{\rm \overline{span}\,}}
\def\cran{{\rm \overline{ran}\,}}
\def\be{\begin{equation}}
\def\ee{\end{equation}}
\def\SZ{Szeg\H{o}}

\begin{document}
\title[ Contractions with rank one defect operators]
{ Contractions with rank one defect operators and truncated CMV
matrices}
\author{
Yury~Arlinski\u{i}}
\address{Department of Mathematical Analysis \\
East Ukrainian National University \\
Kvartal Molodyozhny 20-A \\
Lugansk 91034 \\
Ukraine} \email{yma@snu.edu.ua}

\author{
Leonid~Golinski\u{i}}
\address{Matematics Division, Institute for Low Temperature Physics and
Engineering \\
47 Lenin Ave, Kharkov, 61103} \email{golinskii@ilt.kharkov.ua}

\author{
Eduard~Tsekanovski\u{i}}
\address{Department of Mathematics, P.O.Box 2044 \\
Niagara University, NY 14109, USA } \email{tsekanov@niagara.edu}


\begin{abstract}
The main issue we address in the present paper are the new models
for completely nonunitary contractions with rank one defect
operators acting on some Hilbert space of dimension $N\leq\infty$.
These models complement nicely the well-known models of
Liv$\rm{\check{s}}$ic and Sz.-Nagy--Foias. We show that each such
operator acting on some finite-dimensional (respectively separable
infinite-dimensional Hilbert space) is unitarily equivalent to some
finite (respectively semi-infinite) truncated CMV matrix obtained
from the ``full'' CMV matrix by deleting the first row and the first
column, and acting in $\dC^N$ (respectively $\ell^2(\dN)$). This
result can be viewed as a nonunitary version of the famous
characterization of unitary operators with a simple spectrum due to
Cantero, Moral and Vel\'azquez, as well as an analog for contraction
operators of the result from \cite{ArlTsek} concerning dissipative
non-self-adjoint operators with a rank one imaginary part. It is
shown that another functional model for contractions with rank one
defect operators takes the form of the compression $f(\zeta)\to
P_\cK \left(\zeta\, f(\zeta)\right)$ on the Hilbert space
$L^2(\dT,d\mu)$ with a probability measure $\mu$ onto the subspace
$\cK=L^2(\dT,d\mu)\ominus \dC$. The relationship between
characteristic functions of sub-matrices of the truncated CMV matrix
with rank one defect operators and the corresponding Schur iterates
is established. We develop direct and inverse spectral analysis for
finite and semi-infinite truncated CMV matrices. In particular, we
study the problem of reconstruction of such matrices from their
spectrum or the mixed spectral data involving Schur parameters. It
is pointed out that if the mixed spectral data contains zero
eigenvalue, then no solution, unique solution or infinitely many
solutions may occur in the inverse problem for truncated CMV
matrices. The uniqueness theorem for recovered truncated CMV matrix
from the given mixed spectral data is established. In this part the
paper is closely related to the results of Hochstadt and
Gesztesy--Simon obtained for finite self-adjoint Jacobi matrices.

\end{abstract}
\maketitle
 \tableofcontents
\section{Introduction}

It is well known \cite{AkhGl} that every self-adjoint or unitary
operator with a simple spectrum acting on some separable Hilbert
space is unitarily equivalent to the operator of multiplication by
the independent variable on the Hilbert space $L^2(\dR,d\mu)$ or
$L^2(\dT,d\mu)$, respectively, where $d\mu$ is a probability measure
on the real line $\dR$ or on the unit circle
$\dT=\{\zeta\in\dC:|\zeta|=1\}$. The matrix representation of
self-adjoint operators with simple spectrum was established for the
first time by Stone \cite{Akh}. He proved that every self-adjoint
operator with a simple spectrum is unitarily equivalent to a certain
Jacobi (tri-diagonal) matrix of the form
\begin{equation}
\label{01}
  J=\begin{pmatrix} b_1 & a_1 & 0 &0   & 0 &
\cdot &
\cdot  \\
a_1 & b_2 & a_2 & 0 &0& \cdot &
\cdot  \\
0    & a_2 & b_3 & a_3 &0& \cdot &
\cdot   \\
\cdot & \cdot & \cdot & \cdot & \cdot & \cdot & \cdot
\end{pmatrix},
\end{equation}
where $a_k>0$, and $b_k$ are real numbers for all $k\in\dN$. The
non-self-adjoint version of the Stone theorem has been recently
obtained in \cite{ArlTsek} for dissipative non-self-adjoint
operators with rank one imaginary part. It turned out that the
matrix representation of such operators is a non-self-adjoint Jacobi
matrix of the form (\ref{01}) with only nonreal first entry $b_1$
satisfying $\IM b_1>0$.

The problem of the canonical matrix representation  of a unitary
operator with a simple spectrum has been recently solved by M.
Cantero, L. Moral and L. Vel\'azquez in \cite{CMV1}. They introduced
and studied five-diagonal unitary matrices of the form
\begin{equation}\label{02}
%\[
\cC=\cC(\{\alpha_n\})=\begin{pmatrix}
\bar{\alpha}_0&\bar{\alpha}_1\rho_0&\rho_1\rho_0&0&0&\ldots\cr
\rho_0&-\bar{\alpha}_1\alpha_0&-\rho_1\alpha_0&0&0&\ldots\cr
0&\bar{\alpha}_2\rho_1&-\bar{\alpha}_2\alpha_1&\bar{\alpha}_3\rho_2&\rho_3
\rho_2& \ldots\cr
0&\rho_2\rho_1&-\rho_2\alpha_1&-\bar{\alpha}_3\alpha_2&-\rho_3\alpha_2&\ldots
\cr 0&0&0&\bar{\alpha}_4\rho_3&-\bar{\alpha}_4\alpha_3&\ldots\cr
\ldots&\ldots&\ldots&\ldots&\ldots&\ldots
\end{pmatrix}.
\end{equation}
%\]
Such matrix appears as a matrix representation of the unitary
operator $(Uf)(\zeta)=\zeta f(\zeta)$ in $L_2(\dT,d\mu)$ with
respect to the orthonormal system $\{\chi_n\}$ obtained by
orthonormalization  of the sequence
$\{1,\zeta,\zeta^{-1},\zeta^2,\zeta^{-2},\ldots\}$. The so called
Schur parameters or Verblunsky coefficients $\{\alpha_n\}$,
$|\alpha_n|< 1$, arise in the Szeg\H{o} recurrence formula
$$\zeta \Phi_n(\zeta)=\Phi_{n+1}(\zeta)+
\bar\alpha_n\zeta^n\ovl{\Phi_n(1/\bar\zeta)},\qquad n=0,1, \ldots
$$
for monic orthogonal with respect to $d\mu$ polynomials
$\{\Phi_n\}$, and $\rho_n:=\sqrt{1-|\alpha_n|^2}$. The matrices
$\cC(\{\alpha_n\})$ are called the {\it CMV matrices}. The spectral
analysis of CMV matrices has recently attracted much attention, and
we refer on this matter to the papers \cite{CMV1, CMV2, GZ, S1, S2,
S3}.

The spectral theory of non-self-adjoint and nonunitary operators and
their models is based on the concept of {\it characteristic
function} of the corresponding operator or the operator colligation
\cite{BrR2, Br, Br1, L1, L2, L, LYan, N1, NKh, Pavl, Pol, SF}.

In this paper we employ the Sz.-Nagy--Foias theory \cite{SF} and the
unitary colligations approach \cite{Br1} to the spectral analysis of
contractions acting on Hilbert spaces. The corresponding
characteristic function belongs to the Schur class of
operator-valued functions holomorphic in the open unit disk $\dD$.
By Sz.-Nagy--Foias theorem \cite[Proposition VI.2.1]{SF} each
completely nonunitary contraction $T$ with rank one defect operators
$D_T=(I-T^*T)^{1/2}$ and $D_{T^*}=(I-TT^*)^{1/2}$ (shortly, with
rank one {\it defects}) is unitarily equivalent to the operator
(functional model) of the form
\[
\begin{split}
&\sH_\Theta=\left(H^2\oplus{\rm{clos}}\,\Delta
L^2(\dT))\right)\ominus
\left\{\Theta u\oplus\Delta u:u\in H^2\right\}=\\
&=\left\{\begin{pmatrix}f\cr g \end{pmatrix}: f\in H^2,\;g\in
{\rm{clos}}\,\Delta L^2(\dT)), P_{H^2}(\ovl \Theta f+\Delta g)=0\right\},\\
&\sT_\Theta \begin{pmatrix}f\cr g \end{pmatrix}=P_{\sH_\Theta}\,
\zeta\begin{pmatrix}f\cr g \end{pmatrix},\;\sT^*_\Theta
\begin{pmatrix}f\cr g
\end{pmatrix}=\begin{pmatrix}\bar \zeta (f-f(0))\cr \bar \zeta g
\end{pmatrix} \begin{pmatrix}f\cr g
\end{pmatrix}\in \sH_\Theta,
\end{split}
\]
where $H^2$ is the Hardy space,
$$\Theta=\Theta_T(z)=\left(-T+zD_{T^*}(I-zT^*)^{-1}D_{T}\right)\uphar\sD_{T}$$
is the characteristic function of $T$, $\Delta^2=1-|\Theta|^2$,
$P_{H^2}$ is the orthogonal projection onto $H^2$ in $L^2(\dT)$, and
$P_{\sH_\Theta}$ is the orthogonal projection onto the model space
$\sH_\Theta$.

{\it We obtain a new functional model that complements the above
mentioned Sz.-Nagy--Foias functional model, and show that every
completely nonunitary contraction $T$ with rank one defects is
unitarily equivalent to the compression $f(\zeta)\to P_\cK
\left(\zeta\, f(\zeta)\right)$ on the Hilbert space $L^2(\dT,d\mu)$
with a probability measure $\mu$ onto subspace
$\cK=L^2(\dT,d\mu)\ominus \dC$.}

We study the so called {\it truncated} CMV matrix $\cT$ obtained
from the ``full'' CMV matrix $\cC=\cC(\{\alpha_n\})$ (\ref{02}) by
deleting the first row and the first column:
\[
\cT=\cT (\{\alpha_n\})=\begin{pmatrix}
-\bar{\alpha}_1\alpha_0&-\rho_1\alpha_0&0&0&\ldots\cr
\bar{\alpha}_2\rho_1&-\bar{\alpha}_2\alpha_1&\bar{\alpha}_3\rho_2&\rho_3\rho_2&
\ldots\cr
\rho_2\rho_1&-\rho_2\alpha_1&-\bar{\alpha}_3\alpha_2&-\rho_3\alpha_2&\ldots\cr
%0&0&\bar{\alpha}_4\rho_3&-\bar{\alpha_4}\alpha_3&\ldots\cr
\ldots&\ldots&\ldots&\ldots&\ldots
\end{pmatrix}.\]
In the semi-infinite case $\cT$ takes on the block-matrix form (see
Section 4.3)
\[
  \cT=\begin{pmatrix} \cB_1 & \cC_1 & 0 &0   & 0 &
\cdot &
\cdot  \\
\cA_1 & \cB_2 & \cC_2 & 0 &0& \cdot &
\cdot  \\
0    & \cA_2 & \cB_3 & \cC_3 &0& \cdot &
\cdot   \\
\cdot & \cdot & \cdot & \cdot & \cdot & \cdot & \cdot
\end{pmatrix}.
\]
{\it It turned out that the truncated CMV matrix $\cT(\{\alpha_n\})$
is a contraction with rank one defects, and the Sz.-Nagy--Foias
characteristic function that agrees with the Schur function $f$
associated with the measure $\mu$ \cite{S2}. Moreover, we show that
the sub-matrix $\cT^{(k)}(\{\alpha_n\})$ obtained from
$\cT(\{\alpha_n\})$ by deleting the first $k$ rows and columns is
also a contraction with rank one defects, and its characteristic
function agrees with the well known $k^{th}$ Schur iterate}
\[
f_{k}(z)=\frac{f_{k-1}(z)-\alpha_{k-1}}{z(1-\bar\alpha_{k-1}
f_{k-1}(z))},\qquad f_0(z)=f(z).
\]
This relation is an analog of the corresponding relation between the
$m$-function of a Jacobi matrix and the m-function of its sub-matrix
(cf. \cite{GS}).

{\it Our main result states that each completely nonunitary
contraction $T$ with rank one defects is unitarily equivalent to the
one-parameter family $\cT(\{e^{it}\alpha_n\})$, where $\{\alpha_n\}$
are the Schur parameters of the Sz.-Nagy--Foias characteristic
function of $T$.
 We develop direct and inverse spectral analysis for
finite and  semi-infinite truncated CMV matrices. It is shown that
given an arbitrary set of $N$ not necessarily distinct numbers from
$\dD$ there is a one-parameter family of unitarily equivalent
$N\times N$ truncated CMV matrices having those numbers as the
eigenvalues counting algebraic multiplicity. We prove the uniqueness
of $N\times N$ truncated CMV matrix $\cT$ with given not necessarily
distinct eigenvalues $z_1,\ldots,z_r$ and given  {\it first} $N-r+1$
Schur parameters $\alpha_0(\cT),\ldots,\alpha_{N-r}(\cT)$. This
result on inverse spectral analysis of finite truncated CMV matrices
is an analog of the Hochstadt \cite{H3} and Gesztesy-Simon \cite{GS}
uniqueness theorems for finite self-adjoint Jacobi matrices as well
as for established in \cite{ArlTsek} uniqueness theorem for finite
non-self-adjoint Jacobi matrices with rank one imaginary part. We
obtain the existence of $N\times N$ truncated CMV matrix $\cT$ when
its eigenvalues $z_1,\ldots,z_m$ and the {\it last} Schur parameters
$\alpha_m(\cT),\ldots, \alpha_N(\cT)$ are known.}

Here is a summary of the rest of the paper. In Sections 2 and 3 we
discuss some basics from the Sz.-Nagy--Foias theory and the unitary
colligations with the focus upon the characteristic function and its
properties. Section 4 provides a brief overview of the theory of
orthogonal polynomials on the unit circle and CMV matrices. The main
results concerning truncated CMV matrices and the models of
completely nonunitary contractions with rank one defects are
presented in Section 5 and 6. The final section 7 deals with the
inverse spectral analysis for truncated CMV matrices.

\section{Contractions, unitary colligations, and their characteristic
functions}
\subsection{Contractions and the Sz.- Nagy -- Foias
characteristic functions}

Let $H$ be a separable Hilbert space with the inner product
$(\cdot,\cdot)$. A bounded linear operator $T$ in $H$ is called a
\textit{contraction} if $\|T\|\le 1$ (for the basic properties of
contractions see \cite [Chapter I]{SF}). If $T$ is a contraction
then the operators
\[
D_T:=(I-T^*T)^{1/2},\qquad D_{T^*}:=(I-TT^*)^{1/2}
\]
are called the \textit{defect operators} of $T$ or shortly {\it
defects}, and the subspaces $\sD_T=\cran D_T,$ $\sD_{T^*}=\cran
D_{T^*}$ the \textit{defect subspaces} of $T$. The dimensions
$\dim\sD_T,$ $\dim\sD_{T^*}$ are known as the \textit{defect
numbers} of $T$. Given a pair of numbers $n,n^*=0,1,\ldots,\infty$
it is easy to construct a contraction with $n=\dim\sD_T$,
$n^*=\dim\sD_{T^*}$. Each contraction $T$ acting on a finite
dimensional Hilbert space has equal defect numbers: $n=n^*$.

The defect operators satisfy the following intertwining relations
\begin{equation}
\label{defect} TD_{T}=D_{T^*}T,\qquad T^*D_{T^*}=D_{T}T^*,
\end{equation}
and the block-operators
\[
\begin{pmatrix}-T^*&D_T\cr D_{T^*} &T \end{pmatrix}:
\begin{pmatrix}\sD_{T^*}\cr H\end{pmatrix}
\to \begin{pmatrix}\sD_{T}\cr
H\end{pmatrix},\qquad\begin{pmatrix}-T&D_{T^*}\cr D_{T} &T^*
\end{pmatrix}:\begin{pmatrix}\sD_{T}\cr H\end{pmatrix} \to
\begin{pmatrix}\sD_{T^*}\cr H\end{pmatrix}
\]
are unitary operators in the corresponding orthogonal sums of the
spaces. It follows from \eqref{defect} that
$T\sD_T\subset\sD_{T^*}$, $T^*\sD_{T^*}\subset\sD_T$, and $T(\ker
D_T)=\ker D_{T^*},$ $T^*(\ker D_{T^*})=\ker D_{T}$. Moreover,
$T\uphar\ker D_{T}$ and $T^*\uphar\ker D_{T^*}$ are isometric
operators. It follows that $T$ is a \textit{quasi-unitary extension}
\cite{L1} of the isometric operator $V=T\uphar\ker D_{T}$ (for the
definition see Section 6.2).

A contraction $T$ is called \textit{completely nonunitary} if there
is no nontrivial reducing subspace of $T$, on which $T$ generates a
unitary operator. One of the fundamental results of the contractions
theory \cite [Theorem I.3.2]{SF} reads that, given a contraction $T$
in $H$, there is a canonical orthogonal decomposition
\[
H=H_0\oplus H_1, \qquad T=T_0\oplus T_1, \quad T_j=T\uphar H_j,
\quad j=0,1,
\]
where $H_0$ and $H_1$ reduce $T$, $T_0$ is a completely nonunitary
contraction, and $T_1$ is a unitary operator. Moreover,
\[
H_1= \left(\bigcap\limits_{n\ge 1}\ker
D_{T^n}\right)\bigcap\left(\bigcap\limits_{n\ge 1}\ker
D_{T^{*n}}\right),
\]
so,
\begin{equation}
\label{cu} T \;\mbox{is completely
nonunitary}\;\iff\left(\bigcap\limits_{n\ge 1}\ker
D_{T^n}\right)\bigcap\left(\bigcap\limits_{n\ge 1}\ker
D_{T^{*n}}\right)=\{0\}.
\end{equation}
Clearly,
\begin{equation}
\label{perp}
\begin{split}
&\bigcap\limits_{n\ge 1}\ker D_{T^n}=H\ominus\cspan\left\{T^{*n}D_T
H,\;
n=0,1,\ldots\right\},\\
&\bigcap\limits_{n\ge 1}\ker
D_{T^{*n}}=H\ominus\cspan\left\{T^nD_{T^*} H,\;
n=0,1,\ldots\right\}.
\end{split}
\end{equation}
Let $V$ be an isometry in $H$. A subspace $\Omega$ in $H$ is called
wandering for $V$ if $V^p\Omega\perp V^q\Omega$ for all $p,q\in
\dZ_+$, $p\ne q$. Since $V$ is an isometry, the latter is equivalent
to $V^n\Omega\perp \Omega$ for all $n\in\dN$. If
$H=\sum_{n=0}^\infty\oplus V^n\Omega$, then $V$ is called a
\textit{unilateral shift} and $\Omega$ is called the generating
subspace. The dimension of $\Omega$ is called the multiplicity of
the unilateral shift $V$. It is well known \cite[Theorem I.1.1]{SF}
that $V$ is a unilateral shift if and only if $\bigcap_{n=0}^\infty
V^n H=\{0\}.$ Clearly, if an isometry $V$ is the unilateral shift in
$H$, then $\Omega=H\ominus VH$ is the generating subspace for $V$.

Given a contraction $T$ in $H$ and a subspace $\sH\subset H$, the
unilateral shift $V:\sH\to\sH$ is said to be \textit{contained in
$T$}, if $\sH$ is invariant for $T$, and $T\uphar\sH= V$ \cite{D0}.
 The subspaces $\bigcap\limits_{n\ge 1}\ker D_{T^n}$ and
$\bigcap\limits_{n\ge 1}\ker D_{T^{*n}}$ are invariant for $T$ and
$T^*$, respectively, and the operators $V_T:=T\uphar
\bigcap\limits_{n\ge 1}\ker D_{T^n}$ and
$V_{T^*}:=\uphar\bigcap\limits_{n\ge 1}\ker D_{T^{*n}}$ are
unilateral shifts. Moreover, $V_T$ and $V_{T^*}$ are the maximal
unilateral shifts contained in $T$ and $T^*$. The multiplicities of
the shifts $V_T$ and $V_{T^*}$ do not exceed the defect numbers
$\dim\sD_{T^*}$ and $\dim\sD_{T}$, respectively \cite{DR}. If $T$ is
a completely nonunitary contraction with rank one defects, then (see
\cite{D0}, \cite[Theorem 1.7]{DR})
\begin{equation}
\label{shift}
\begin{split}
& T\;\mbox{does not contain the unilateral shift}\iff
T^*\;\mbox{does
not contain the unilateral shift}\\
& \iff\bigcap\limits_{n\ge 1}\ker D_{T^n}=\{0\}\iff
\bigcap\limits_{n\ge 1}\ker D_{T^{*n}}=\{0\}.
\end{split}
\end{equation}

 The function (see \cite[Chapter VI]{SF})
\[
\Theta_T(z)=\left(-T+zD_{T^*}(I-zT^*)^{-1}D_{T}\right)\uphar\sD_{T}
\]
is known as the \textit{characteristic function} of the Sz.-Nagy --
Foias type of a contraction $T$. This function belongs to the
\textit{Schur class} ${\bf S}(\sD_T,\sD_{T^*})$ of
$\cL(\sD_T,\sD_{T^*})$-valued holomorphic in the unit disk $\dD$
operator-functions, i.e., $\|\Theta_T(z)\|\le 1$ for $z\in\dD$.
Moreover, the function $\Theta_T$ satisfies the condition
$\|\Theta_T(0)f\|<\|f\|$ for all $f\in\sD_T\setminus\{0\}$. The
characteristic functions of $T$ and $T^*$ are connected by the
relation
%\begin{equation}
\[
%\label{conj}
 \Theta_{T^*}(z)=\Theta^*_T(\bar z),\quad z\in\dD.
%\end{equation}
\]

Two operator-valued functions $\Theta_1\in{\bf S}(\sM_1,\sN_1)$ and
$\Theta_2\in{\bf S}(\sM_2,\sN_2)$ are said to \textit{agree} if
there are two unitary operators $V:\sN_1\to\sN_2$ and
$W:\sM_2\to\sM_1$ such that
\[
V\Theta_1(z)W=\Theta_2(z),\quad z\in\dD.
\]
It is well known \cite[Theorem VI.3.4]{SF}, that two completely
nonunitary contractions $T_1$ and $T_2$ are unitary equivalent if
and only if their characteristic functions $\Theta_{T_1}$ and
$\Theta_{T_2}$ agree.

Every operator-valued function $\Theta$ from the Schur class ${\bf
S}(\sM,\sN)$ has almost everywhere nontangential strong limit values
$\Theta(\zeta)$, $\zeta \in \dT$. A function $\Theta \in {\bf
S}(\sM,\sN)$ is called \textit{inner} if
$\Theta^*(\zeta)\Theta(\zeta)=I_\sM$ for a.e. $\zeta\in\dT$, and
\textit{co-inner} if $\Theta(\zeta)\Theta^*(\zeta)=I_\sN$ for a.e.
$\zeta\in\dT$. A function $\Theta\in {\bf S}(\sM,\sN)$ is called
\textit{bi-inner}, if it is both inner and co-inner. A contraction
$T$ on a Hilbert space $\sH$ belongs to the classes $C_{0\,\cdot}$
($C_{\cdot\, 0}$), if
\[
s-\lim\limits_{n\to\infty}T^n=0  \qquad
(s-\lim\limits_{n\to\infty}T^{*n}=0),
\]
respectively. By definition $C_{00}:=C_{0\,\cdot}\cap C_{\cdot\,
0}$. The completely nonunitary part of a contraction $T$ belongs to
the class $C_{\cdot\,0}$, $C_{0\,\cdot}$, or $C_{00}$  if and only
if its characteristic function $\Theta_T(z)$  is inner, co-inner, or
bi-inner, respectively (cf. \cite[Section VI.2]{SF}).

In the following statement \cite[Theorem VI.4.1]{SF} the spectrum of
completely nonunitary contractions is described.
\begin{theorem}
\label{Spectr} Let $T$ be a completely nonunitary contraction on
$H$. Denote by $S_T$ the set of points $z \in \dD$ for which the
operator $\Theta_T(z)$ is not boundedly invertible, together with
those $z\in \dT$ not lying on any of the open arcs of $\dT$ on which
$\Theta_T$ is a unitary operator valued analytic function.
Furthermore, denote by $S_T^0$ the set of points $z\in \dD$ for
which $\Theta_T(z)$ is not invertible at all. Then the spectrum
$\sigma(T)$ of $T$ agrees with $S_T$, and the point spectrum
$\sigma_p(T)$ with $S_T^0$.
\end{theorem}

If $T$ is a completely nonunitary contraction with rank one defects,
and if $z_0$ is an eigenvalue of $T$, then the geometric
multiplicity of $z_0$ is one, the algebraic multiplicity is finite,
and the characteristic function $\Theta_T$ admits the following
factorization
\[
\Theta_T(z)=c\prod\frac{\bar z_k}{z_k}\,\frac{z_k-z}{1-\bar z_k
z}\,\exp\left(-\int\limits_0^{2\pi}\frac{e^{it}+z}{e^{it}-z}d\mu(t)\right)\,\
\exp\left(\frac{1}{2\pi}
\int\limits_0^{2\pi}\frac{e^{it}+z}{e^{it}-z}\ln k(t)dt\right),
\]
where $|c|=1$, $k(t)\ge 0$, $\ln k(t)\in L_1[0,2\pi]$, $\mu$ is a
finite nonnegative measure singular with respect to the Lebesgue
measure, and $\{z_k\}$ are the eigenvalues of $T$. In addition, if
$\dim H=N<\infty$, and $T$ is a completely nonunitary contraction in
$H$ with rank one defects, then its characteristic function is the
finite Blaschke product of order $N$ of the form
\[
b(z)=e^{i\f}\prod\limits_{k=1}^m\left(\frac{z-z_k}{1-\bar z_k
z}\right)^{l_k},
\]
where $z_1,\ldots,z_m$ are distinct eigenvalues of $T$ with the
algebraic multiplicities $l_1,\ldots, l_m$, respectively,
$l_1+\ldots+l_m=N$, and $\f\in [0,2\pi)$. Hence, a
finite-dimensional completely nonunitary contraction $T$ with rank
one defects belongs to the class $C_{00}$, and
$\lim\limits_{n\to\infty}||T^n||=0.$ It is easily seen from Theorem
\ref{Spectr} that the point spectrum of a contraction $T$ with rank
one defects agrees with $\dD$ if and only if $\Theta_T\equiv 0$.

\subsection{Unitary colligations and their
characteristic functions} Every contraction $T$ acting on Hilbert
space $H$ can be included into the \textit{unitary operator
colligation} \cite{Br1} \footnote{also known as the
\textit{conservative system \cite{A}}}
\[
\Delta=\left\{\begin{pmatrix} S&G\cr F&T\end{pmatrix};\sM,\sN,
H\right\},
\]
where $\sM$ and $\sN$ are separable Hilbert spaces, and
\[
U=\begin{pmatrix} S&G\cr F&T\end{pmatrix}:\begin{pmatrix}\sM\cr H
\end{pmatrix}\to \begin{pmatrix}\sN\cr H
\end{pmatrix}
\]
is a unitary operator. $T$ is called the \textit{basic operator} of
the unitary colligation $\Delta$. The spaces $\sM$ and $\sN$ are
called the \textit{left outer space} and \textit{right outer space},
respectively. The unitarity of $U$ means
\[
U^*U=\begin{pmatrix}I_\sM&0\cr 0& I_H\end{pmatrix},\qquad
UU^*=\begin{pmatrix}I_\sN&0\cr 0& I_H\end{pmatrix}.
\]
or equivalently,
\begin{equation}
\label{conserv}
\begin{split}
&T^*T+G^*G=I_H,\quad F^*F+S^*S=I_\sM,\quad T^*F+G^*S=0,\\
&TT^*+FF^*=I_H,\quad GG^*+SS^*=I_\sN,\quad TG^*+FS^*=0.
\end{split}
\end{equation}
The colligation
\begin{equation}
\label{delta0} \Delta_0=\left\{\begin{pmatrix}-T^*&D_T\cr D_{T^*} &T
\end{pmatrix};\;\sD_{T},\sD_{T^*}, H\right\}
\end{equation}
provides an example of the unitary colligation with given basic
operator $T$.

Let $\Delta=\left\{\begin{pmatrix} S&G\cr F&T\end{pmatrix};
\sM,\sN,H\right\}$ be a unitary colligation. Define the following
subspaces in $H$
\begin{equation}
\label{co} H^{(c)}=\cspan\{T^nF\sM,\; n=0,1,\ldots\},\qquad
H^{(o)}=\cspan\{T^{*n}G^*\sN,\; n=0,1,\ldots\}.
\end{equation}
The subspaces $H^{(c)}$ and $H^{(o)}$ are called the
\textit{controllable} and the \textit{observable} subspaces,
respectively. Let
\begin{equation}
\label{ORTCOMP}
 (H^{(c)})^\perp:=H\ominus
H^{(c)},\qquad (H^{(o)})^\perp :=H\ominus H^{(o)}.
\end{equation}
 A unitary colligation $\Delta$ is called
\textit{prime} if $\overline{H^{(c)}+H^{(o)}}=H$. Clearly, the
latter condition is equivalent to
\[
(H^{(c)})^\perp\cap(H^{(o)})^\perp=\{0\},
\]
 From \eqref{conserv} and \eqref{ORTCOMP} we get
\begin{equation}
\label{SP}
\begin{split}
&(H^{(c)})^\perp=\bigcap\limits_{n\ge 0}\ker
(F^*T^{*n})=\bigcap\limits_{n\ge
0}\ker(D_{T^*}T^{*n})=\bigcap\limits_{n\ge 1}\ker(D_{T^{*n}}),\\
& (H^{(o)})^\perp=\bigcap\limits_{n\ge 0}\ker
(GT^n)=\bigcap\limits_{n\ge 0}\ker(D_{T}T^n)=\bigcap\limits_{n\ge
1}\ker(D_{T^n}).
\end{split}
\end{equation}
It follows now from \eqref{cu} that \textit{the unitary colligation
$\Delta=\left\{\begin{pmatrix} S&G\cr F&T\end{pmatrix};
\sM,\sN,H\right\}$ is prime if and only if $T$ is a completely
nonunitary operator}.

Given a unitary colligation $\Delta=\left\{\begin{pmatrix} S&G\cr
F&T\end{pmatrix}; \sM,\sN,H\right\}$, its characteristic function
\footnote{\textit{the transfer function} of the system \cite{A}}
\cite[Section 3]{Br1} is defined by
\[
\Theta_\Delta(z)=S+zG(I_H-zT)^{-1}F,\quad z\in\dD.
\]
This function belongs to the Schur class ${\bf S}(\sM,\sN)$ of
$\cL(\sM,\sN)$-valued holomorphic in the unit disk $\dD$
operator-functions. In particular, the characteristic function of
the unitary colligation $\Delta_0$ \eqref{delta0}
\[
\Theta_0(z)=\left(-T^*+zD_{T}(I-zT)^{-1}D_{T^*}\right)\uphar\sD_{T^*}
\]
is in fact the Sz.-Nagy -- Foias characteristic function of the
operator $T^*$.

Two prime unitary colligations
\[
\Delta_1=\left\{\begin{pmatrix} S&G_1\cr F_1&T_1\end{pmatrix};
\sM,\sN,H_1\right\}\quad\mbox{and}\quad
\Delta_2=\left\{\begin{pmatrix} S&G_2\cr F_2&T_2\end{pmatrix};
\sM,\sN,H_2\right\}
\]
which have equal characteristic functions are unitarily equivalent
in the following sense \cite[Theorem 3.2]{Br1}: there exists a
unitary operator $V:H_1\to H_2$ such that
\[
VT_1=T_2V, \;VF_1=F_2,\; G_2V=G_1\iff \begin{pmatrix} I_\sN&0\cr 0
&V\end{pmatrix}\begin{pmatrix} S&G_1\cr
F_1&T_1\end{pmatrix}=\begin{pmatrix} S&G_2\cr
F_2&T_2\end{pmatrix}\begin{pmatrix} I_\sM&0\cr 0 &V\end{pmatrix}
\]
Besides, given $\Theta\in {\bf S}(\sM,\sN)$, there exists a prime
unitary colligation
\[
\Delta=\left\{\begin{pmatrix} S&G\cr F&T\end{pmatrix};
\sM,\sN,H\right\}
\]
such that $\Theta_\Delta=\Theta$ in $\dD$ \cite[Theorem 5.1]{Br1}.

Later on in Section 3 we will need the following result.
\begin{theorem}
\label{CH1} Let $T$ be a contraction acting on Hilbert space $H$
with finite defect numbers. Suppose that $\sM$ and $\sN$ are two
given Hilbert spaces such that $\dim \sN=\dim\sD_T$ and
$\dim\sM=\dim\sD_{T^*}$. Then all unitary colligations with the
basic operator $T$ and outer subspaces $\sM$ and $\sN$ take the form
\begin{equation}
\label{UN}\Delta=\left\{\begin{pmatrix} -KT^*M&KD_T\cr
D_{T^*}M&T\end{pmatrix}; \sM,\sN,H\right\},
\end{equation}
where $K:\sD_T\to\sN$ and $M:\sM\to\sD_{T^*}$ are unitary operators.
The characteristic function of $\Delta$ is
\[
\Theta_\Delta(z)=K\Theta_{T^*}(z)M,\quad z\in\dD,
\]
i.e., $\Theta_\Delta$ agrees with the characteristic function
$\Theta_{T^*}$ of \ $T^*$.
\end{theorem}

\begin{proof}
Let $\Delta=\left\{\begin{pmatrix} S&G\cr F&T\end{pmatrix};
\sM,\sN,H\right\}$ be a unitary colligation. From the relation
$G^*G+T^*T=I_H$ it follows that
\[
||Gf||^2=||D_{T}f||^2,\quad f\in H.
\]
Hence, the operator $K:\sD_T\to \sN$ defined by
\[
KD_Tf=Gf,\quad f\in H,
\]
is isometric, and $\ran K=\sN$. Similarly, the relation
$FF^*+TT^*=I_H$ yields that the operator $N:\sD_{T^*}\to\sM$ given
by the relation
\[
ND_{T^*}f=F^*f,\quad f\in H
\]
is isometric, and $\ran N=\sM$. So  $M=N^*:\sM\to\sD_{T^*}$ is
unitary, and $F=D_{T^*}M$.

From the relation $T^*F+G^*S=0$ we get $T^*D_{T^*}M+D_{T}K^*S=0$.
Hence by \eqref{defect} $T^*M+K^*S=0$. As $\ran M=\sD_{T^*}$, $\ran
K^*=\sD_T$, and $T\sD_{T^*}\subset\sD_T$, we have
\[
S=-KT^*M.
\]
Observe also that
\[
\begin{split}
&TG^*+FS^*=TD_{T}K^*-D_{T^*}MM^*TK^*=0\\
&SS^*+GG^*=KT^*MM^*TK^*+KD^2_{T}K^*=K(T^*T+I-T^*T)K^*=I_\sN,\\
&S^*S+F^*F=M^*TK^*KT^*M+M^*D_{T^*}M=M^*(TT^*+I-TT^*)M=I_\sM.
\end{split}
\]
Thus, all conditions \eqref{conserv} are satisfied, i.e., the
colligation $\Delta$ is of the form \eqref{UN}.

Conversely, if $\dim\sN=\dim\sD_T<\infty$,
$\dim\sM=\dim\sD_{T^*}<\infty$, and $K:\sD_T\to\sN$ and
$M:\sM\to\sD_{T^*}$ are unitary operators, then one can easily see
that $U=\begin{pmatrix} -KT^*M&KD_T\cr
D_{T^*}M&T\end{pmatrix}:\begin{pmatrix}\sM\cr H\end{pmatrix}\to
\begin{pmatrix}\sN\cr H\end{pmatrix}$ is a unitary operator, i.e.,
the relations \eqref{conserv} are satisfied. It follows that
$$\Delta=\left\{\begin{pmatrix} S&G\cr
F&T\end{pmatrix}; \sM,\sN,H\right\}$$ is a unitary colligation,
where $G=KD_T,$ $F=D_{T^*}M$, $S=-KT^*M$.

For the characteristic function $\Theta_\Delta$ we obtain for all
$z\in\dD$
\[
\Theta_\Delta(z)=S+zG(I-zT)^{-1}F=-KT^*M+zKD_{T}(I-zT)^{-1}D_{T^*}M=K\Theta_
{T^*}(z)M.
\]
\end{proof}
\begin{corollary}
\label{COL1} Let $T$ be a contraction with finite defect numbers,
$\dim \sN=\dim\sD_T$, $\dim\sM=\dim\sD_{T^*}$, and let
$\Delta=\left\{\begin{pmatrix} S&G\cr F&T\end{pmatrix};
\sM,\sN,H\right\}$ be a unitary colligation. Then all other unitary
colligations with the basic operator $T$ and outer subspaces $\sM$
and $\sN$ take the form
\[
\wt \Delta=\left\{\begin{pmatrix} C_1SC_2&C_1G\cr
FC_2&T\end{pmatrix}; \sM,\sN,H\right\},
\]
where $C_1$ and $C_2$ are unitary operators in $\sN$ and $\sM$,
respectively.
\end{corollary}
\begin{proof} By Theorem \ref{CH1} we have
\[
G=KD_T,\;F=D_{T^*}M,\; S=-KT^*M,
\]
where $K:\sD_T\to\sM$ and $M:\sN\to\sD_{T^*}$ are unitary operators.
If $\wt\Delta=\left\{\begin{pmatrix} \wt S&\wt G\cr \wt
F&T\end{pmatrix}; \sM,\sN,H\right\}$ is some other unitary
colligation then $\wt G=\wt KD_T,\;\wt F=D_{T^*}\wt M,\; \wt S=-\wt
KT^*\wt M,$  where $\wt K:\sD_T\to\sM$ and $\wt M:\sN\to\sD_{T^*}$
are unitary operators. Let $C_1:=\wt KK^{-1},\; C_2:=M^{-1}\wt M$.
Then $C_1$ and $C_2$ are unitary operators in $\sN$ and $\sM$,
respectively, and
\[
\wt G=C_1 G,\quad \wt F=FC_2,\quad \wt S=C_1SC_2,\] as needed.
\end{proof}

\section{Completely nonunitary contractions with rank one defects
and the corresponding unitary colligations}

\begin{theorem}
\label{cyclic} Each contraction $T$ with rank one defects on the
Hilbert space $H$ can be included into the unitary colligation
\[
\Delta=\left\{\begin{pmatrix} S&G\cr F&T\end{pmatrix};\ \dC,\dC,
H\right\}.
\]
Let $\vec 1=\begin{pmatrix}1\cr 0\end{pmatrix}\in\dC\oplus H$, and
let the subspaces $(H^{(c)})^\perp$ and $(H^{(o)})^\perp$ in $H$ be
defined by \eqref{ORTCOMP}.
 Then
\begin{equation}
\label{coperp} \begin{split} & (H^{(c)})^\perp=(\dC\oplus
H)\ominus\cspan\{U^n\vec 1;\; n=0,1,\ldots\},\\
& (H^{(o)})^\perp=(\dC\oplus H)\ominus\cspan\{U^{*n}\vec 1;\;
n=0,1,\ldots\},
\end{split}
\end{equation}
and so the following conditions are equivalent:
\begin{enumerate}
\def\labelenumi{\rm (\roman{enumi})}
\item the unitary colligation $\Delta=\left\{\begin{pmatrix}
S&G\cr F&T\end{pmatrix};\ \dC,\dC, H\right\}$ is prime; \item $T$ is
a completely nonunitary contraction; \item $\vec 1$ is the cyclic
vector for $U$: $\cspan\{U^n \vec 1,\; n\in\dZ\}=\dC\oplus H$.
\end{enumerate}
All other unitary colligations with the basic operator $T$ and the
outer spaces $\dC$ are of the form
\begin{equation}
\label{OTHER} \wt \Delta=\left\{\begin{pmatrix} c_1c_2S&c_1G\cr
c_2F&T\end{pmatrix};\ \dC,\dC,H\right\},
\end{equation}
where $|c_1|=|c_2|=1$.
\end{theorem}
\begin{proof} Since $\dim\sD_T=\dim\sD_{T^*}=1$, by Theorem \ref{CH1}
we can choose the unitary colligation $\Delta=\left\{\begin{pmatrix}
S&G\cr F&T\end{pmatrix};\dC,\dC, H\right\}$ of the form \eqref{UN},
i.e., $S=-KT^*M ,\;G=KD_T,\; F=D_{T^*}M$, and $K:\ran D_T\to \dC$,
$M:\dC\to \ran D_{T^*}$ are isometric operators. So
$U=\begin{pmatrix} S&G\cr F&T\end{pmatrix}:\begin{pmatrix}\dC\cr H
\end{pmatrix}\to \begin{pmatrix}\dC\cr H
\end{pmatrix}$ is the unitary operator.

To prove \eqref{coperp}, suppose that the vector $\vec
h=\begin{pmatrix} z\cr h\end{pmatrix}\in\dC\oplus H$ is orthogonal
to the subspace $\cspan\{U^n \vec 1,\; n=0,1,\ldots\}$. Then $
U^{*n}\vec h\perp \vec 1$, $n=0,1,\ldots$ so $z=0$ and $\vec
h=\begin{pmatrix}0 \cr h\end{pmatrix}$. By using
$U^*=\begin{pmatrix} S^*&F^*\cr G^*&T^*\end{pmatrix}$, we get
consequently
\[
F^*h=0,\quad F^*T^*h=0,\quad F^*T^{*2}h=0,\quad\ldots, \quad
F^*T^{*k}h=0,\ldots
\]
It follows from (\ref{SP}) that $h\in (H^{(c)})^\perp.$ Conversely,
if $h \in (H^{(c)})^\perp$ then $h\perp \cspan\{U^n \vec 1,\;
n=0,1,\ldots\}$. Similarly, $(H^{(o)})^\perp=(\dC\oplus
H)\ominus\left(\cspan\{U^{*n}\vec 1,\; n=0,1,\ldots\}\right),$  as
needed.

We arrive at the following conclusion:
\[
\begin{split}
&\vec 1\:\mbox{is a cyclic vector for}\; U \iff (H^{(c)})^\perp\cap
(H^{(o)})^\perp=\{0\}\iff\\
&\mbox{the unitary colligation}\; \Delta=\left\{\begin{pmatrix}
S&G\cr
F&T\end{pmatrix};\dC,\dC, H\right\}\;\mbox{is prime}\iff \\
&\mbox{the operator}\; T\;\mbox{is completely nonunitary}.
\end{split}
\]
By Corollary \ref{COL1} all other unitary colligations with the
basic operator $T$ and the outer subspace $\dC$ are given by
\eqref{OTHER} with $|c_1|=|c_2|=1.$
\end{proof}
\begin{remark}
\label{VD1} In terms of the Naimark dilations of a probability
operator-valued measure on the unit circle, the main result of
Theorem \ref{cyclic} is proved in \cite[Theorem 1.20]{D1}.
\end{remark}
Let us give more precise expressions for the operators $F,G$, and
$S$. Let $\wh\f_1\in\sD_T$, $\wh\f_2\in \sD_{T^*}$. Put
\[
\f_1=\frac{\wh\f_1}{\|\wh\f_1\|},\qquad\f_2=\frac{\wh\f_2}{\|\wh\f_2\|}.
\]
Then
\[
\begin{split}
&Kh=b_1(h,\f_1),\quad h\in\ran D_T,\\
&M^*g=b_2(g,\f_2),\quad g\in\ran D_{T^*},
\end{split}
\]
where $|b_1|=|b_2|=1$. Observe that $T\f_1=-\alpha_0 \f_2$ and
$T^*\f_2=-\bar{\alpha}_0\f_1$, where $\alpha_0$ is a complex number
from $\dD$. It follows that
\[
D^2_T\f_1=(1-|\alpha_0|^2)\f_1,\qquad
D^2_{T^*}\f_2=(1-|\alpha_0|^2)\f_2.
\]
Let $\rho_0=\sqrt{1-|\alpha_0|^2}.$  Since $\dim(\ran
D^2_T)=\dim(\ran D^2_{T^*})=1$, the number $\rho_0$ is a unique
positive eigenvalue of $D_T $($D_{T^*}$).
%From \eqref{GN} it follows that
Next,
\[
\begin{split}
&Gh=b_1(D_T h,\f_1)=b_1(h, D_T\f_1)=b_1\rho_0(h,\f_1),\\
&F^*h=b_2(D_{T^*}h,\f_2)=b_2(h, D_{T^*}\f_2)=b_2\rho_0(h,\f_2),\quad
h\in H.
\end{split}
\]
Hence $F 1=\rho_0\bar b_2\f_2$. Since $S=-KT^*M$, we get
\[
S1=-b_1 \bar b_2(T^*\f_2,\f_1)=b_1 \bar b_2\bar\alpha_0.
\]
In the case $\dim H=N<\infty$ the operator $T$ can be given by the
$N\times N$ matrix with respect to some orthonormal basis and  we
can choose $\wh\f_1$(respectively, $\wh\f_2)$, as one of the nonzero
columns of the matrix $I-T^*T\;$  $(I-TT^*)$. In addition,
\[
Trace(I-T^*T)=Trace(I-TT^*)=\rho^2_0.
\]
Thus, if
$\f_2=\begin{pmatrix}\f^{(1)}_{2}\cr\f^{(2)}_{2}\cr\ldots\cr
\f^{(N)}_{2}
   \end{pmatrix}$, then the column $F$ takes the form
$F=\bar b_2\rho_0\begin{pmatrix}\f^{(1)}_{2}
\cr\f^{(2)}_{2}\cr\ldots\cr \f^{(N)}_{2} \end{pmatrix}$. \newline If
$\f_1=\begin{pmatrix}\f^{(1)}_{1}\cr\f^{(2)}_{1}\cr\ldots\cr
\f^{(N)}_{1}
   \end{pmatrix}$, then
   the row $G$ takes the form  $
G=b_1\rho_0
\begin{pmatrix}\bar\f^{(1)}_{1}&\bar\f^{(2)}_{1}&\ldots&\bar\f^{(N)}_{1}
\end{pmatrix}.
$ Finally, the number $S$ is given by $-b_1 \bar b_2(T^*\f_2,\f_1)$.

If $\dim H=N$ and $T$ is a completely nonunitary contraction with
rank one defects, then $\Theta_\Delta$ is a finite Blaschke product
\[
\Theta_\Delta(z)=e^{i\f}\prod\limits_{k=1}^N\frac{z-\bar z_k}{1-
z_k\, z},
\]
where the numbers $ z_1,\dots,  z_N$ are the eigenvalues of $T$.
Since all other unitary colligations are of the form \eqref{OTHER},
for the characteristic function $\Theta_{\wt\Delta}(z)$ we get
$\Theta_{\wt\Delta}(z)=c_1c_2\Theta_\Delta(z)=e^{it}\Theta_\Delta(z),$
$z\in\dD$, and $t\in[0,2\pi).$

Let $U$ be a unitary operator with a cyclic vector $e$, acting on
the Hilbert space $H$. The spectral measure $\mu$ associated with
$U$ and $e$ provides the relation
\[
\left(F(U)e,e\right)=\int_{\dT} F(\zeta)d\mu(\zeta),
\]
which is the Spectral Theorem for unitaries. For instance,
\be\label{cfun} F(z)=\left((U+zI)(U-zI)^{-1}e,e\right)=
\int_{\dT}\frac{\zeta+z}{\zeta-z}\,d\mu(\zeta)\,, \quad z\in\dD \ee
is the Carath\'eodory function (\ref{carat}), i.e., $F$ is
holomorphic in the unit disc $\dD$, $\RE F>0$ in $\dD$, and
$F(0)=1$.
\begin{theorem}
\label{char}  Let $T$ be a completely nonunitary contraction with
rank one defects, $\Delta=\left\{\begin{pmatrix} S&G\cr
F&T\end{pmatrix};\ \dC,\dC,H\right\}$ be the prime unitary
colligation, and $\Theta_\Delta$ be its characteristic function. Put
\begin{equation}
\label{Carat}
 F(z)=\left((U+zI)(U-zI)^{-1}\vec 1,\vec
1\right),\quad z\in\dD,
\end{equation}
where $U=\begin{pmatrix} S&G\cr
F&T\end{pmatrix}:\begin{pmatrix}\dC\cr H
\end{pmatrix}\to \begin{pmatrix}\dC\cr H
\end{pmatrix}$.
Then
\begin{equation}
\label {REL} \ovl{\Theta_\Delta(\bar
z)}=\frac{1}{z}\,\frac{F(z)-1}{F(z)+1},\qquad
F(z)=\frac{1+z\ovl{\Theta_\Delta(\bar
z)}}{1-z\ovl{\Theta_\Delta(\bar z)}},\quad z\in\dD.
\end{equation}
%In other words, $\ovl{\Theta_\Delta(\bar z)}=f(z)$, where the left
%hand side is the characteristic function of $T$, and the right hand
%side is the Schur function for the spectral measure $d\mu$
%associated with $U$ and $\vec 1$.

%Therefore, the entries of the corresponding truncated CMV matrix
%$\cT$ unitary equivalent to $T$ are determined by the Schur
%parameters $\{\alpha_n\}$ of the function $\ovl{\Theta_\Delta(\bar
%z)}$ whose agrees with characteristic function of the operator
%$T$. All other truncated CMV matrices unitary equivalent to $T$
%have parameters  $\{e^{i\f}\alpha_k\}$, where $\f\in [0,2\pi).$
\end{theorem}
\begin{proof} We use the well known Schur--Frobenius formula for the inverse
of block operators (see, e.g., \cite[Section 0.2]{D}, \cite[p.
57]{Gant}). Let $\sH_1$ and $\sH_2$ be two Hilbert spaces, and
$\Phi$ an operator in $\sH_1\oplus\sH_2$ given by the block operator
matrix
\[
\Phi=\begin{pmatrix}A& B\cr C&D \end{pmatrix}: \begin{pmatrix}
\sH_1\cr\sH_2\end{pmatrix}\to \begin{pmatrix}
\sH_1\cr\sH_2\end{pmatrix}.
\]
Suppose that $D^{-1}\in\cL(\sH_2)$ and
$(A-BD^{-1}C)^{-1}\in\cL(\sH_1)$. Then $\Phi^{-1}\in
\cL(\sH_1\oplus\sH_2,\sH_1\oplus\sH_2)$ and
\[
\Phi^{-1}=\begin{pmatrix}K^{-1}& -K^{-1}BD^{-1}\cr
-D^{-1}CK^{-1}&D^{-1}+D^{-1}CK^{-1}BD^{-1}
\end{pmatrix},
\]
where $K=A-BD^{-1}C$.

Applying this formula for
\[
\Phi=I-zU=\begin{pmatrix}1-zS&-zG\cr-zF&I-zT\end{pmatrix}:\begin{pmatrix}
\dC\cr H\end{pmatrix}\to \begin{pmatrix} \dC\cr H\end{pmatrix},
\quad z\in\dD,
\]
we get $K=1-zS-z^2G(I-zT)^{-1}F=1-z\Theta_\Delta(z)$. Therefore
\[
\left((I-zU)^{-1}\vec 1,\vec 1\right)=\frac{1}{1-z\Theta_\Delta(z)},
\quad z\in\dD.
\]
Let
\[
\Psi(z)=\left((I+z U)(I-zU)^{-1}\vec 1,\vec 1\right),\quad z\in\dD.
\]
Clearly, the equality  $F(z)=\ovl{\Psi(\bar z)}$ holds, which yields
\eqref{REL}.
\end{proof}
\begin{remark}
\label{VD2} Relations \eqref{REL} is proved in \cite[Theorem 1.20,
Comments 2.8]{D1}. Our proof is different.
\end{remark}

\section{OPUC and CMV matrices}
\subsection{Basics of OPUC}
It is well recognized now that the theory of orthogonal polynomials
on the real line plays an important role in the spectral theory of
self-adjoint operators (and close to such operators) acting on
Hilbert spaces. Likewise, the theory of orthogonal polynomials on
the unit circle (OPUC) appears in the same fashion in the study of
unitary operators and close to such operators. Here we recall some
rudiments and advances of the OPUC theory.

If $\mu$ is a nontrivial probability measure on $\dT$ (that is, not
supported on a finite set), the monic orthogonal polynomials
$\Phi_n(z,\mu)$ (or $\Phi_n$ if $\mu$ is understood) are uniquely
determined by \be\label{monic} \Phi_n(z)=\prod_{j=1}^n (z-z_{n,j}),
\qquad \int_{\dT} \zeta^{-j}\Phi_n(\zeta)\,d\mu=0, \quad
j=0,1,\ldots,n-1, \ee so on the Hilbert space $L^2(\dT, d\mu)$,
$\langle\Phi_n, \Phi_m\rangle=0$, $n\not=m$. We also consider the
orthonormal polynomials $\phi_n$ of the form
$\phi_n=\Phi_n/\|\Phi_n\|$.

In case when $\mu$ is supported on a finite set, that is,
\be\label{triv} \mu=\sum_{k=1}^N \mu_k\delta(\zeta_k), \quad
\zeta_k\in\dT, \ee a finite number of orthogonal polynomials
$\{\Phi_k\}_{k=0}^{N-1}$ can be defined in the same manner.

Clearly, (\ref{monic}) and the fact that the space of polynomials of
degree at most $n$ has dimension $n+1$ imply \be\label{monic1}
\deg(P)=n, \quad P\bot\,\zeta^j, \quad j=0,1,\ldots, n-1 \Rightarrow
P=c\Phi_n. \ee

On $L^2(\dT, d\mu)$ the anti-unitary map
$f^*(\zeta):=\zeta^n\ovl{f(\zeta)}$ (which depends on $n$) is
naturally defined. The set of polynomials of degree at most $n$ is
left invariant: \be\label{starpoly} P(z)=\sum_{j=0}^n p_jz^j
\Rightarrow P^*(z)=\sum_{j=0}^n \bar p_{n-j}z^j. \ee (\ref{monic1})
now implies \be\label{starort} \deg(P)\leq n, \quad P\bot\,\zeta^j,
\quad j=1,\ldots, n \Rightarrow P=c\Phi_n^*. \ee

A key feature of the unit circle is that the multiplication $Uf=z f$
in $L^2(\dT, d\mu)$ is a unitary operator. So the difference
$\Phi_{n+1}(z)-z\Phi_n(z)$ is of degree $n$ and orthogonal to $z^j$
for $j=1,2,\ldots,n$, and by (\ref{starort}) \be\label{szego}
\Phi_{n+1}(z)=z\Phi_n(z)-\bar\alpha_n(\mu)\Phi_n^*(z) \ee with some
complex numbers $\alpha_n(\mu)$, called the {\it Verblunsky
coefficients} \cite{S2}. (\ref{szego}) is known as the {\it \SZ\
recurrences} after its first occurrence in the celebrated book
\cite{Sz} of G. Szeg\H{o}. (\ref{szego}) at $z=0$ imply
\be\label{vc} \alpha_n(\mu)=\alpha_n=-\ovl{\Phi_{n+1}(0)}. \ee It is
known that for nontrivial measures $|\alpha_n|<1$ for all
$n=0,1,2,\ldots$, and for trivial measures (\ref{triv}) one has a
finite set of Verblunsky coefficients $\{\alpha_n\}_{n=0}^{N-1}$
with $|\alpha_n|<1$, $n=0,1,\ldots,N-2$, and $|\alpha_{N-1}|=1$.
Since it arises often, define \be\label{rho}
\rho_j:=\sqrt{1-|\alpha_j|^2}, \qquad 0\leq\rho_j\leq 1, \qquad
|\alpha_j|^2+\rho_j^2=1. \ee The inverse \SZ\ recurrences are also
of interest (cf. \cite [Theorem 1.5.4]{S2}): \be\label{invszego}
z\Phi_n(z)=\rho_n^{-2}\left(\Phi_{n+1}(z)+\bar\alpha_n\Phi_{n+1}^*(z)\right)
. \ee The norm of the polynomials $\Phi_n$ in $L^2(\dT,d\mu)$ can be
computed by:
\[
||\Phi_n||=\prod\limits_{j=0}^{n-1}\rho_j,\qquad n=1,2,\ldots.
\]

Let $\dD^\infty$ be the set of complex sequences
$\{\alpha_j\}_{j=0}^\infty$ with $|\alpha_j|<1$. The map $\cS$, from
$\mu\rightarrow\{\alpha_j(\mu)\}_{j=0}^\infty$, is a well defined
map from the set $\cP$ of nontrivial probability measures on $\dT$
to $\dD^\infty$. It was S.~Verblunsky who proved that $\cS$ is a
bijection. As a matter of fact, $\cS$ is a homeomorphism, provided
$\cP$ is equipped with the weak*-topology, and $\dD^\infty$ with the
topology of component convergence. Moreover, it follows directly
from (\ref{szego}) that for two measures $\mu_1$ and $\mu_2$
\[
\alpha_j(\mu_1)=\alpha_j(\mu_2), \quad j=0,1,\ldots,n-1 \
\Rightarrow \ \Phi_j(z,\mu_1)=\Phi_j(z,\mu_2), \quad j=0,1,\ldots,n.
\]
Conversely, by (\ref{invszego})
\[
\Phi_n(z,\mu_1)=\Phi_n(z,\mu_2) \ \Rightarrow \
\alpha_j(\mu_1)=\alpha_j(\mu_2), \quad j=0,1,\ldots,n-1.
\]

The orthonormal set $\{\phi_n\}_{n\geq0}$ does not necessarily form
a basis in $L^2(\dT, d\mu)$ (e.g., if $d\mu=dm$ is the normalized
Lebesgue measure on $\dT$, then $\phi_n=\zeta^n$ and $\zeta^{-1}$ is
orthogonal to all $\phi_n$). A celebrated result of Szeg\H{o} --
Kolmogorov -- Krein reads that $\{\phi_n\}$ is a basis in $L^2(\dT,
d\mu)$ if and only if $\log \mu'\not\in L^1(\dT)$, where $\mu'$ is
the Radon -- Nikodym derivative of $\mu$ with respect to $dm$. In
addition, the following result holds true (cf. \cite[Theorem
1.5.7]{S2}).
\begin{theorem}
\label{complete} For any nontrivial probability measure $\mu$ on the
unit circle, the following are equivalent
\begin{enumerate}
\def\labelenumi{\rm (\roman{enumi})}
\item $\lim\limits_{n\to\infty}||\Phi_n||=0;$ \item
$\sum\limits_{n=0}^\infty |\alpha_n|^2=\infty;$ \item the system
$\{\phi_n\}_{n=0}^\infty$ is the orthonormal basis in $L^2(\dT,
d\mu)$.
\end{enumerate}
\end{theorem}
Note that if $\sum\limits_{n=0}^\infty |\alpha_n|^2<\infty$  and $P$
is the orthogonal projection in $L^2(\dT, d\mu)$ onto
$\cspan\{\zeta^n,\; n=0,1,\ldots\}$, then (see \cite{S1})
\begin{equation}
\label{PROJ}
 \|(I-P)\bar \zeta\|=\prod\limits_{n=0}^\infty\rho_n.
\end{equation}

Let us now turn to the basic properties of zeros
$\{z_{n,j}\}_{j=1}^n$ of OPUC. It is well known (cf., e.g., \cite
[Theorem 1.7.1]{S2}) that $|z_{n,j}|<1$ for all $n$ and $j$.
Moreover, a result of Geronimus \cite [Theorem 1.7.5]{S2} reads that
given a monic polynomial $P_n$ of degree $n$ with all its zeros
inside $\dD$, there is a (nontrivial) probability measure $\mu$ on
$\dT$ such that $P_n=\Phi_n(\mu)$. Actually, there are infinitely
many such measures, all of them have the same Verblunsky
coefficients up to the order $n-1$, and the same moments up to the
order $n$. Given a monic polynomial $P_n$ with all its zeros inside
the disk, let us call a monic polynomial $Q_{n+m}$ an {\it
extension} of $P_n$, if there is a measure $\mu$ such that
\[
P_n=\Phi_n(\mu), \quad Q_{n+m}=\Phi_{n+m}(\mu).
\]
To obtain all such extensions one just has to extend a sequence of
Verblunsky coefficients $\alpha_0,\ldots,\alpha_{n-1}$, which are
completely determined by $P_n$, by a sequence
$\beta_0,\ldots,\beta_{m-1}$ with arbitrary $\beta_j\in\dD$ and then
apply (\ref{szego}).

One of the most recent advances in the study of zeros of OPUC is the
theorem of Simon and Totik \cite [Theorem 1.7.15]{S2}, which claims
that given a polynomial $P_n$ as above, and an arbitrary set of
points $z_1,\ldots,z_m$ in the unit disk, not necessarily distinct,
there is an extension $Q_{n+m}$ of $P_n$ such that $Q_{n+m}(z_j)=0$,
$j=1,2,\ldots,m$, counting the multiplicity. The latter as usual
means that
\[
z_k=z_{k+1}=\ldots=z_{k+p}\Rightarrow
Q_{n+m}(z_k)=Q_{n+m}'(z_k)=\ldots=Q_{n+m}^{(p)}(z_k)=0.
\]
The uniqueness of such extension is an open problem. A particular
case $m=1$ appeared earlier in \cite{AV}. Now $\beta_0=\alpha_n$ is
defined uniquely from (\ref{szego}) by
\[
0=Q_{n+1}(z_1)=z_1P_n(z_1)-\bar\alpha_nP_n^*(z_1).
\]
This result will play a key role in the inverse problems with mixed
data in Section 7.

\subsection{Geronimus theory}
There is an important analytic aspect of the OPUC theory which was
developed by Geronimus \cite{ger1, ger2} in 1940's.

Given a probability measure $\mu$ on $\dT$, define the {\it
Carath\'eodory function} by \be\label{carat} F(z)=
F(z,\mu):=\int_{\dT}\frac{\zeta+z}{\zeta-z}\,d\mu(\zeta)=
1+2\sum_{n=1}^\infty \beta_nz^n, \quad \beta_n=\int_{\dT}
\zeta^{-n}d\mu \ee the moments of $\mu$. $F$ is an analytic function
in $\dD$ which obeys $\RE F>0$, $F(0)=1$. The {\it Schur function}
is then defined by \be\label{schur} f(z)=
f(z,\mu):=\frac{1}{z}\,\frac{F(z)-1}{F(z)+1}, \qquad F(z)=\frac{1+z
f(z)}{1-z f(z)}\,, \ee so it is an analytic function in $\dD$ with
$\sup_{\dD}|f(z)|\leq1$. A one-to-one correspondence can be easily
set up between the three classes (probability measures,
Carath\'eodory and Schur functions). Under this correspondence $\mu$
is trivial, that is, supported on a finite set, if and only if the
associate Schur function is a finite Blaschke product. Moreover,
this Blaschke product has the order $N-1$ for measures (\ref{triv}).

We proceed with the Schur algorithm. Given a Schur function $f=f_0$,
which is not a finite Blaschke product, define inductively
  \be\label{schuralg}
  f_{n+1}(z)=\frac{f_n(z)-\gamma_n}{z(1-\bar\gamma_n f_n(z))}\,, \qquad
  \gamma_n=f_n(0). \ee
It is clear that the sequence $\{f_n\}$ is an {\it infinite}
sequence of Schur functions (called the $n^{th}$ Schur iterates) and
neither of its terms is a finite Blaschke product. The numbers
$\{\gamma_n\}$ are called the {\it Schur parameters:}
\[
\cS f=\{\gamma_0,\gamma_1,\ldots,\}.
\]
In case when
\[
f(z)=e^{i\f}\prod_{k=1}^N \frac{z-z_k}{1-\bar z_k z}
\]
is a finite Blaschke product of order $N$, the Schur algorithm
terminates at the N-th step. The sequence of Schur parameters
$\{\gamma_k\}_{k=0}^N$ is finite, $|\gamma_k|<1$ for
$k=0,1,\ldots,N-1$, and $|\gamma_N|=1$.

If a Schur function $f$ is not a finite Blaschke product, the
connection between the non-tangential limit values $f(\zeta)$
 and its Schur parameters $\{\gamma_n\}$ is given by the
formula
\begin{equation}
\label{L1} \prod\limits_{n=0}^\infty
(1-|\gamma_n|^2)=\exp\left\{\int_\dT\ln(1-|f(\zeta)|^2)dm\right\}
\end{equation}
(see \cite {Boy}). It follows that
\[
\sum\limits_{n=0}^\infty|\gamma_n|^2=\infty\iff
\ln(1-|f(\zeta)|^2)\notin L^1(\dT).
\]
In addition, if one of the conditions
%\begin{equation}
%\label{innerfu}
\begin{enumerate}
\item $\limsup_{n\to\infty}|\gamma_n|=1,$ \item
 $\lim_{n\to\infty}\gamma_n\gamma_{n+m}=0$ for each
$m=1,2,\ldots$, but $\limsup_{n\to\infty}|\gamma_n|>0$,
\end{enumerate}
%\end{equation}
is fulfilled, then $f$ is the inner function (see \cite{Rakh,
Khfund}).

Later in Section 7 we will make use of the following fundamental
result of Schur \cite{Schur}: the set of all Schur functions $f$
with prescribed first Schur parameters $\gamma_0,\ldots,\gamma_n$ is
given by the linear fractional transformation \be \label{set}
f(z)=\frac{A(z)+zB^*(z)s(z)}{B(z)+zA^*(z)s(z)}\,, \ee where $s$ is
an arbitrary Schur function, and $A,B$ are polynomials of degree at
most $n$. Moreover,
\[
\cS f=\{\gamma_0,\ldots,\gamma_n,\gamma_0(s),\gamma_1(s),\ldots\}
\]
The pair $(A,B)$, known as the Wall pair, is completely determined
by $\{\gamma_j\}_{j=0}^n$. Specifically,
\[
W(z):=\begin{pmatrix}zB^*(z)& A(z)\cr zA^*(z)& B(z)\end{pmatrix}=
Q_{\gamma_0}(z)\,Q_{\gamma_1}(z)\cdots Q_{\gamma_n}(z),
\]
where
\[
Q_\omega(z)=\frac{1}{\sqrt{1-|\omega|^2}}\,\begin{pmatrix}z&\omega\cr
z\bar\omega&\, 1\end{pmatrix}, \quad \omega\in\dD.
\]
By computing determinants, we see that
\[
B^*(z)B(z)-A^*(z)A(z)=z^n\prod_{j=0}^n(1-|\gamma_j|^2)^{1/2},
\]
so $A$ and $B$ have no common zeros in $\dC\setminus\{0\}$. In fact
they have no common zeros at all since $B(0)=1$. It is known also
that $B\not=0$ in $\ovl\dD$, and both $AB^{-1}$ and $A^*B^{-1}$ are
Schur functions.

A straightforward computation shows that $Q_\omega$ (and hence $W$)
are $j$-inner matrix functions:
\[
\begin{split}
&W^*(z)jW(z)\ge j\quad\mbox{for}\quad z\in\dD,\\
&W^*(z)jW(z)=j\quad\mbox{for}\quad z\in\dT
\end{split}
\]
with the signature matrix
\[
j=\begin{pmatrix}-1&0\cr 0& 1 \end{pmatrix}.
\]
For further properties of the Wall pairs see \cite[Section
4]{Khfund}, \cite[Chapter 1.3.8]{S2}.

A curious situation when the Schur parameters for a finite Blaschke
product can be computed explicitly was found by Khrushchev \cite
[formula (1.12)]{Kh}. Let $\mu$ be a nontrivial probability measure
(or measure of the form (\ref{triv}) with big enough $N$) with
Verblunsky coefficients $\{\alpha_k\}$, and $\Phi_n$ be its $n$th
monic orthogonal polynomial. Consider the following Blaschke product
of order $n$
\[
b_0(z):=\frac{\Phi_n(z)}{\Phi_n^*(z)}=\prod_{j=1}^n\frac{z-z_{n,j}}{1-\bar
z_{n,j} z}, \qquad b_0(0)=-\bar\alpha_{n-1}.
\]
 It is a matter of a simple computation based on (\ref{invszego})
to make sure that
\[
b_{1}(z)=\frac{b_0(z)-b_0(0)}{z(1-\bar b_0(0)
b_0(z))}=\frac{\Phi_{n-1}(z)}{\Phi_{n-1}^*(z)}\,.
\]
Hence the Schur parameters of $b_0$ are of the form \be\label{kh}
\cS
b_0=\{-\bar\alpha_{n-1},-\bar\alpha_{n-2},\ldots,-\bar\alpha_0,1\}.
\ee

The fundamental paper of Schur \cite{Schur} had appeared a few years
before Szeg\H{o} introduced the notion of orthogonal polynomials on
the unit circle. Amazingly, neither of them benefited from the ideas
of the other. Only 20 years later Geronimus put them together and
came up with the following fundamental result (see \cite [Theorem
IX, p. 111]{ger1})

\begin{theorem}
\label{TT} Let $\mu$ be a nontrivial probability measure on $\dT$
and $f$ its Schur function with the Schur parameters $\gamma_n(f)$.
Then $\gamma_n(f)=\alpha_n(\mu)$. For measures (\ref{triv}) the
latter equality holds for $n=0,1,\ldots,N-1$.
\end{theorem}
It is clear now why a minus and conjugate is taken in (\ref{szego}).

We complete with the result which will be used later on in Section
7.
\begin{theorem}
\label{TTT} Given two sets $\alpha_0,\ldots,\alpha_{n-1}$ and
$z_1,\ldots,z_m$ of complex numbers in $\dD$, and $\gamma\in\dT$,
there exists a finite Blaschke product $b$ of order $n+m$ such that
\begin{enumerate}
\def\labelenumi{\rm (\roman{enumi})}
    \item $\cS
    b=\{\omega_0,\ldots,\omega_{m-1},\alpha_0,\ldots,\alpha_{n-1},\gamma\}$,
    \item $b(z_j)=0$, $j=1,\ldots,m$, counting multiplicity.
\end{enumerate}
\end{theorem}
\begin{proof}
Denote $\beta_k:=-\bar\gamma\bar\alpha_{n-k-1}$, $k=0,1,\ldots,n-1$
and construct a system of monic orthogonal polynomials
$\{\Phi_k(z,\beta)\}_{k=0}^n$  by (\ref{szego}). The theorem of
Simon-Totik claims that there is a measure $\mu$ with
\[
\Phi_n(z,\mu)=\Phi_n(z,\beta), \qquad \Phi_{n+m}(z_j,\mu)=0, \quad
j=1,\ldots,m,
\]
counting the multiplicity. The first equality means that
$\alpha_k(\mu)=\beta_k$, $k=1,\ldots,n-1$. Finally, put
\[
b(z):=\gamma\,\frac{\Phi_{n+m}(z,\mu)}{\Phi_{n+m}^*(z,\mu)}\,.
\]
The result now follows from Khrushchev's formula (\ref{kh}).
\end{proof}

Note that for $m=1$ the Blaschke product is uniquely determined.

\subsection{CMV matrices}

One of the most interesting developments in the OPUC theory in
recent years is the discovery by Cantero, Moral, and Vel\'azquez
\cite{CMV1, CMV2} of a matrix realization for the operator of
multiplication by $\zeta$ on $L^2(\dT,d\mu)$ which is a unitary
matrix of finite band size (i.e.,
$|\langle\zeta\chi_m,\chi_n\rangle|=0$ if $|m-n|>k$ for some $k$);
in this case, $k=2$ to be compared with $k=1$ for the Jacobi
matrices, which correspond to the real line case. The CMV basis
(complete, orthonormal system) $\{\chi_n\}$ is obtained by
orthonormalizing the sequence
$1,\zeta,\zeta^{-1},\zeta^2,\zeta^{-2},\ldots$, and the matrix,
called the {\it CMV matrix},
$$ \cC={\cC}(\mu)=\|c_{n,m}\|_{m,n=0}^\infty =\|\langle \zeta\chi_m,\chi_n
\rangle\|,
    \qquad m,n\in\dZ_+  $$
is five-diagonal. Remarkably, the $\chi$'s can be expressed in terms
of $\phi$'s and $\phi^*$'s:
$$ \chi_{2n}(z)=z^{-n}\phi_{2n}^*(z), \qquad
\chi_{2n+1}(z)=z^{-n}\phi_{2n+1}(z), \quad n\in\dZ_+, $$ and the
matrix elements in terms of $\alpha$'s and $\rho$'s:
\begin{equation}
\label{cmvmatr} \cC=\cC(\{\alpha_n\})=\begin{pmatrix}
\bar{\alpha}_0&\bar{\alpha}_1\rho_0&\rho_1\rho_0&0&0&\ldots\cr
\rho_0&-\bar{\alpha}_1\alpha_0&-\rho_1\alpha_0&0&0&\ldots\cr
0&\bar{\alpha}_2\rho_1&-\bar{\alpha}_2\alpha_1&\bar{\alpha}_3\rho_2&\rho_3
\rho_2& \ldots\cr
0&\rho_2\rho_1&-\rho_2\alpha_1&-\bar{\alpha}_3\alpha_2&-\rho_3\alpha_2&\ldots
\cr 0&0&0&\bar{\alpha}_4\rho_3&-\bar{\alpha}_4\alpha_3&\ldots\cr
\ldots&\ldots&\ldots&\ldots&\ldots&\ldots
\end{pmatrix},
\end{equation}
$\alpha$'s are the Verblunsky coefficients and $\rho$'s are given in
(\ref{rho}).

It is not hard to write down a general formula for the matrix
entries $c_{ij}$ (see \cite{Gol}). Let $2\epsilon_m:=1-(-1)^m$,\
$m\in\mathbb{Z}_+$, and $\epsilon_{-1}=1$, so $\{\epsilon_m\}_{m\geq
0}=\{0,1,0,1,\ldots\}$,
$$ \epsilon_m+\epsilon_{m+1}=1, \qquad \epsilon_m\epsilon_{m+1}=0,
   \qquad \epsilon_m-\epsilon_{m+1}=(-1)^{m+1}. $$
Then
\begin{equation}\label{cmventr1} \begin{split}
c_{mm} &= -\ovl\alpha_m\alpha_{m-1},\\
c_{m+2,m} &= \rho_m\rho_{m+1}\epsilon_m, \qquad
c_{m,m+2}=\rho_m\rho_{m+1}\epsilon_{m+1},
\end{split}\end{equation}
and
\begin{equation}\label{cmventr2} \begin{split}
  c_{m+1,m} &=
\ovl\alpha_{m+1}\rho_m\epsilon_m-\alpha_{m-1}\rho_m\epsilon_{m+1},\\
  c_{m,m+1} &=
\ovl\alpha_{m+1}\rho_m\epsilon_{m+1}-\alpha_{m-1}\rho_m\epsilon_m.
\end{split}\end{equation}
It is clear (cf. \cite[Theorem 1]{BD}), that any semi-infinite CMV
matrix $\cC$ (\ref{cmvmatr}) can be written in the three-diagonal
block-matrix form

\begin{equation}
\label{THREE}
  \cC=\begin{pmatrix} \cB_0 & \cC_0 & 0 &0   & 0 &
\cdot &
\cdot  \\
\cA_0 & \cB_1 & \cC_1 & 0 &0& \cdot &
\cdot  \\
0    & \cA_1 & \cB_2 & \cC_2 &0& \cdot &
\cdot   \\
\cdot & \cdot & \cdot & \cdot & \cdot & \cdot & \cdot
\end{pmatrix}
\end{equation}
with
\begin{equation}
\label{BLOCKS} \begin{split} \cB_0
&=\begin{pmatrix}\bar\alpha_0\end{pmatrix},\quad
\cC_0=\begin{pmatrix}\bar\alpha_1\rho_0&\rho_1\rho_0\end{pmatrix},\quad
\cA_0=\begin{pmatrix}\rho_0\\0\end{pmatrix},\\
\cA_n
&=\begin{pmatrix}\rho_{2n}\rho_{2n-1}&-\rho_{2n}\alpha_{2n-1}\cr 0&0
\end{pmatrix},\quad
\cB_n=\begin{pmatrix}-\bar\alpha_{2n-1}\alpha_{2n-2}&-\rho_{2n-1}\alpha_{2n-2}\cr
\bar\alpha_{2n}\rho_{2n-1}&-\bar\alpha_{2n}\alpha_{2n-1}
\end{pmatrix},\\
\cC_n
&=\begin{pmatrix}0&0\cr\bar\alpha_{2n+1}\rho_{2n}&\rho_{2n+1}\rho_{2n}
\end{pmatrix}, \qquad n=1,2,\ldots.
\end{split}
\end{equation}
There is a nice multiplicative structure of the CMV matrices. In the
semi-infinite case $\cC$ is the product of two matrices:
$\cC=\cL\cM$, where
\begin{equation}
\label{lm}
\begin{split}
\cL
&=\Psi(\alpha_0)\oplus\Psi(\alpha_2)\oplus\ldots\Psi(\alpha_{2m})
\oplus\ldots,\\
\cM &={\bf 1}_{1\times 1}\oplus\Psi(\alpha_1)\oplus\Psi(\alpha_3)
\oplus\ldots\oplus\Psi(\alpha_{2m+1})\oplus\ldots,
\end{split}
\end{equation}
and
$\Psi(\alpha)=\begin{pmatrix}\bar\alpha&\rho\cr\rho&-\alpha\end{pmatrix}.$
The finite $(N+1)\times (N+1)$ CMV matrix $\cC$ obeys
$\alpha_0,\alpha_1,\ldots,\alpha_{N-1}\in \dD$ and $ |\alpha_{N}|=1$
is also the product $\cC=\cL\cM$, where in this case
$\Psi(\alpha_{N})=\left(\bar\alpha_{N}\right)$.

It is just natural to take the ordered set
$1,\zeta^{-1},\zeta,\zeta^{-2},\zeta^2,\ldots$ instead of
$1,\zeta,\zeta^{-1},\zeta^2,\zeta^{-2},\ldots$, that leads to the
alternate CMV basis $\{x_n\}$ and the alternate CMV matrix
\begin{equation}
\label{altcmvmatr} \wt\cC=\|\langle \zeta x_m,x_n
\rangle\|=\begin{pmatrix} \bar{\alpha}_0&\rho_0&0&0&0&\ldots\cr
\bar{\alpha_1}\rho_0&-\bar{\alpha}_1\alpha_0&\bar{\alpha}_2\rho_1&\rho_2\rho_1&0&\ldots\cr
\rho_1\rho_0&-\rho_1\alpha_0&-\bar{\alpha}_2\alpha_1&-\rho_2\alpha_1&0&\ldots\cr
0&0&\bar{\alpha}_3\rho_2&-\bar{\alpha}_3\alpha_2&\bar{\alpha_4}\rho_3&\ldots
\cr
0&0&\rho_3\rho_2&-\rho_3\alpha_2&-\bar{\alpha_4}\alpha_3&\ldots\cr
\ldots&\ldots&\ldots&\ldots&\ldots&\ldots
\end{pmatrix},
\end{equation}
which turns out to be the transpose of $\cC$ (see \cite[Corollary
4.2.6]{S2}). Furthermore, $\cL=\cL^t$ and $\cM=\cM^t$ imply
$\wt\cC=\cC^t=\cM\cL$.

An important relation between CMV matrices and monic orthogonal
polynomials similar to the well-known property of orthogonal
polynomials on the real line
%\begin{equation}\label{20.8}
\[
\Phi_n(z)=\det(zI_n-{\cC}^{(n)})
\]
%\end{equation}
holds, where ${\cC}^{(n)}$ is the principal $n\times n$ block of
${\cC}$.

One of the most important results of Cantero, Moral, and Vel\'azquez
\cite{CMV1} states that each unitary operator $U$ with the simple
spectrum (i.e. having a cyclic vector $e_1$) acting on some
infinite-dimensional separable Hilbert space (respectively,
finite-dimensional Hilbert space) is unitarily equivalent to a
certain CMV matrix in $l_2(\dZ_+)$ (respectively, in $\dC^n$). The
corresponding $\alpha$'s come up as the Verblunsky coefficients of
the spectral measure $d\mu$ of $U$ associated with $e_1$. This is
the analog of Stone's self-adjoint cyclic model theorem. To be more
precise, let us, following \cite{S3}, call a {\it cyclic unitary
model} a unitary operator $U$ acting on a separable Hilbert space
$\cH$ with the distinguished cyclic unit vector $v_0$. Two cyclic
unitary models, $(\cH,U,v_0)$ and $(\tilde\cH, \tilde U, \tilde
v_0)$ are called equivalent if there is a unitary operator $W$ from
$\cH$ onto $\tilde\cH$ such that $Wv_0=\tilde v_0$ and
$WUW^{-1}=\tilde U$. It is clear that $\delta_0=(1,0,0,\ldots)^t$ is
cyclic for any CMV matrix $\cC$. Moreover, every class of equivalent
unitary models contains exactly one CMV model
$(\ell^2,\cC,\delta_0)$.

\section{A model in the space $L^2(\dT, d\mu)$ of a completely nonunitary contraction
with rank one defects} \label{s5}
\begin{theorem}
\label{L2} Let $T$ be a completely nonunitary contraction with rank
one defects. Then there exists a probability measure $\mu$ on $\dT$
such that $T$ is unitarily equivalent to the following operator
\begin{equation} \label{UNMODL2}
\sT h(\zeta)=P_\sH\left(\zeta h(\zeta)\right),\qquad
h\in\sH:=L^2(\dT,d\mu)\ominus\dC,
\end{equation}
where $P_\sH$ is the orthogonal projection in $L^2(\dT,d\mu)$ onto
$\sH$. The Schur function associated with $\mu$ is exactly the
characteristic function of $T$.
\end{theorem}
\begin{proof}
Include $T$ into a prime unitary colligation
$\Delta=\left\{\begin{pmatrix} S&G\cr F&T\end{pmatrix};
\dC,\dC,H\right\}$. The characteristic function $\Theta_\Delta$
agrees with the characteristic function of $T^*$. By Theorem
\ref{cyclic} the vector $\vec 1=\begin{pmatrix}1\cr 0\end{pmatrix}$
is cyclic for the unitary operator $U=\begin{pmatrix} S&G\cr
F&T\end{pmatrix}$.

Let $E_U(\zeta)$ be the resolution of identity for $U$. Define
$d\mu(\zeta):=(dE_U(\zeta)\vec 1,\vec 1)$ and put
\[
\cU f(\zeta)=\zeta f(\zeta)
 \]
the unitary multiplication operator in $L^2(\dT, d\mu)$. By the
spectral theorem for unitaries with cyclic vectors (cf.
\cite[Section 1.4.5]{S2}) there exists a unitary operator
$W:\dC\oplus H\to L^2(\dT,d\mu)$ such that
\[
U=W^{-1}\cU W \;\mbox{and}\; W\vec 1=1.
\]
It follows that $W$ takes the block-operator form
\[
W=\begin{pmatrix} 1&0\cr 0&V \end{pmatrix}:\begin{pmatrix}\dC\cr H
\end{pmatrix}\to \begin{pmatrix}\dC\cr\sH
\end{pmatrix},
\]
where $\sH=L^2(\dT,d\mu)\ominus\dC$, $V:H\to
L^2(\dT,d\mu)\ominus\dC$ is a unitary operator. If $\sT$ is given by
\eqref{UNMODL2}, then
\[
\sT:=P_\sH\cU\uphar\sH=VTV^{-1},
\]
i.e., $T$ is unitarily equivalent to $\sT$. Clearly, $\cU$ has the
block form
\[
\cU=\begin{pmatrix}P_\dC \cU\uphar\dC &P_\dC \cU\uphar\sH\cr
P_\sH\cU\uphar\sH&\sT
\end{pmatrix},
\]
where $P_\dC$ is the orthogonal projection in $L^2(\dT,d\mu)$ onto
the subspace $\dC$ of the constant functions in $L^2(\dT,d\mu)$. The
unitary colligation $\Delta$ is unitarily equivalent to the unitary
colligation
\begin{equation}
\label{unmod} \left\{\begin{pmatrix}P_\dC \cU\uphar\dC &P_\dC
\cU\uphar\sH\cr P_\sH\cU\uphar\sH&\sT
\end{pmatrix},\dC,\dC,\sH\right\}.
\end{equation}
Note that
\[
P_\dC(\cU\, 1)=\int_\dT\zeta \,d\mu,\; P_\sH(\cU\,
1)=\zeta-\int_\dT\zeta \,d\mu,\;
P_\dC(\cU^*\,1)=\bar\zeta-\int_\dT\bar\zeta \,d\mu.
\]
 Let
$F(z)=\left((U+zI)(U-zI)^{-1}\vec 1,\vec 1\right)$. Then
\[
F(z)=\left((\cU+zI)(\cU-zI)^{-1}
1,1\right)=\int_{\dT}\frac{\zeta+z}{\zeta-z}\,d\mu(\zeta),
\]
 i.e., $F$ is the Carath\'eodory function associated with $\mu$.
 From Theorem \ref{char} we conclude
\[
\ovl{\Theta_\Delta(\bar z)}=\frac{1}{z}\frac{F(z)-1}{F(z)+1},
\]
and so by \eqref{schur} $\ovl{\Theta_\Delta(\bar z)}$ agrees with
the Schur function associated with $\mu$.
\end{proof}
Let $\{\Phi_n\}$ be the system of monic polynomials orthogonal with
respect to $\mu$, and let $\{\alpha_n\}$ be the corresponding
Verblunsky coefficients. By Geronimus' theorem $\{\alpha_n\}$ are
the Schur parameters of $f$. Let $\sH^{(c)}$ be the controllable
subspace of the unitary colligation \eqref{unmod}. From
\eqref{coperp} it follows that
\[
(\sH^{(c)})^\perp=L^2(\dT, d\mu)\ominus\cspan\{\zeta^n,\;
n=0,1,\ldots\}.
\]
If $\mu$ is a nontrivial measure, then in view of \eqref{PROJ} we
obtain
\[
\|P_{(\sH^{(c)})^\perp}
\bar\zeta\|=\prod\limits_{n=0}^\infty(1-|\alpha_n|^2)^{1/2}.
\]
The latter is equivalent to
\[
\|P_{(\sH^{(c)})^\perp}P_\dC(\cU^*
\,1)\|=\prod\limits_{n=0}^\infty(1-|\alpha_n|^2)^{1/2}.
\]
Hence, from \eqref{UN} and \eqref{co} we have the equivalence
\begin {equation}
\label{poln}
 \cspan \{\sT^n\sD_{\sT^*},\;n=0,1,\ldots\}=\sH\iff
\sum\limits_{n=0}^\infty|\alpha_n|^2=\infty.
\end{equation}

\begin{remark}
\label{Anmodel} By the construction of Theorem \ref{L2}, the Schur
function $f$ associated with $\mu$ is exactly
$\ovl{\Theta_\Delta(\bar z)}$. Another (unitary equivalent) models
of $T$ are connected with the operators $U_\lambda=\begin{pmatrix}
\bar\lambda S&G\cr \bar\lambda F&T\end{pmatrix}$, where
$|\lambda|=1$. The characteristic function of the unitary
colligation
$$
\Delta_\lambda=\left\{\begin{pmatrix} \bar\lambda S&G\cr \bar\lambda
F&T\end{pmatrix},\dC,\dC, H\right\}
$$
is $\bar\lambda\Theta_\Delta$. The model operator $\sT_\lambda$
takes the form
\[
\sH_\lambda=L^2(\dT,d\mu_\lambda)\ominus\dC,\; \sT_\lambda
h(\zeta)=P_{\sH_\lambda}\left(\zeta h(\zeta)\right),\quad
h(\zeta)\in\sH_\lambda.
\]
The Schur function $f_\lambda$ associated with $\mu_\lambda$ is
$f_\lambda=\lambda\, f$. The connection between the Carath\'eodory
functions $F_\lambda(z)=\left((U_\lambda+zI)(U_\lambda-zI)^{-1}\vec
1,\vec 1\right)$ and $F$ is given by
\[
F_\lambda(z)=\frac{(1-\lambda)+(1+\lambda)
F(z)}{(1+\lambda)+(1-\lambda)F(z)}.
\]
The measures $\mu_\lambda$ are known as the \textit{Aleksandrov
measures} associated with $\mu$ \cite[Section 1.3.9]{S2}.
\end{remark}

\section{Truncated CMV matrices}
\subsection{Truncated CMV matrix as a model for contractions with
rank one defects} \label{ss6}

Let $\cC=\cC(\{\alpha_n\})$ be the CMV matrix given by
\eqref{cmvmatr}. Recall that $\cC(\{\alpha_n\})$ is the matrix
representation of the unitary operator $\cU$ of multiplication by
$\zeta$ in $L^2(\dT, d\mu)$, where $\mu$ is the probability measure
with Verblunsky coefficients $\{\alpha_n\}$. By the Geronimus
theorem the Schur parameters of the Schur function \eqref{schur}
associated with $\mu$ are $\{\alpha_n\}$.

The matrix $\cC$ determines the unitary operator in the space
$l_2(\dZ_+)$ (respectively in $\dC^{N+1}$ in the case of
$(N+1)\times (N+1)$ matrix). The vector $\delta_0=(1,0,0,\ldots)^t$
is cyclic for $\cC$. Consider the matrix
\begin{equation}
\label{MODEL} \cT=\cT (\{\alpha_n\})=\begin{pmatrix}
-\bar{\alpha}_1\alpha_0&-\rho_1\alpha_0&0&0&\ldots\cr
\bar{\alpha}_2\rho_1&-\bar{\alpha}_2\alpha_1&\bar{\alpha}_3\rho_2&\rho_3\rho_2&
\ldots\cr
\rho_2\rho_1&-\rho_2\alpha_1&-\bar{\alpha}_3\alpha_2&-\rho_3\alpha_2&\ldots\cr
0&0&\bar{\alpha}_4\rho_3&-\bar{\alpha}_4\alpha_3&\ldots\cr
\ldots&\ldots&\ldots&\ldots&\ldots
\end{pmatrix}
\end{equation}
obtained from $\cC$ by deleting the first row and the first column.
It is clear from \eqref{THREE} that a semi-infinite $\cT$ takes on
the three-diagonal $2\times 2$ block-matrix form
\[
  \cT=\begin{pmatrix} \cB_1 & \cC_1 & 0 &0   & 0 &
\cdot &
\cdot  \\
\cA_1 & \cB_2 & \cC_2 & 0 &0& \cdot &
\cdot  \\
0    & \cA_2 & \cB_3 & \cC_3 &0& \cdot &
\cdot   \\
\cdot & \cdot & \cdot & \cdot & \cdot & \cdot & \cdot
\end{pmatrix},
\]
where $\cA_n$, $\cB_n$, and $\cC_n$ are defined in \eqref{BLOCKS}.
Henceforth $\cT$ is called a \textit{truncated} CMV matrix. $\cT$ is
the matrix of the operator $\sT=P_\sH\cU\uphar\sH$, where $P_\sH$ is
the orthogonal projection in $L^2(\dT, d\mu)$ onto the subspace
$\sH=L^2(\dT, d\mu)\ominus\dC$.

It is easy to see that given $\cT$ \eqref{MODEL}, the values
$\alpha_n$ are uniquely determined. Indeed, from $(2,2)$ and $(3,2)$
entries we have by (\ref{rho}) $
|\alpha_1|^2=|\bar\alpha_2\alpha_1|^2+\rho_2^2|\alpha_1|^2, $ so
$|\alpha_1|$ and $\rho_1>0$ are known, and we find $\alpha_0$,
$\alpha_2$ from $(1,2)$ and $(2,1)$ entries of (\ref{MODEL}). From
$(2,1)$ and $(2,2)$ entries we get $\rho_2>0$, then $\alpha_1$,
$\alpha_3$, etc. We call $\alpha_n=\alpha_n(\cT)$ the
\textit{parameters} of $\cT$ (\ref{MODEL}).

As it was mentioned in Section 4.3, $\cC=\cL\cM$, $\cL$ and $\cM$
are defined in \eqref{lm}. Given a matrix $A$, we denote by $A_r$
($A_c$) the matrix obtained from $A$ by deleting the first row
(column). Clearly, $A_{rc}={(A_r)}_c$. So we have
$\cT=\cC_{rc}=\cL_r\cM_c$. $\cM_c$ is isometric with $\dim\ran(I-
\cM_c\cM^*_c)=1$, whereas $\cL_r$ is coisometric with
$\dim\ran(I-\cL^*_r\cL_r)=1$.

Let $P_{{\delta_0}^\perp}$ be the orthogonal projection in
$l_2(\dZ_+)$ ($\dC^{N+1})$ onto the subspace ${\delta_0}^\perp\cong
l_2(\dN)$ $(\dC^{N})$. Then the matrix $\cT$ determines on the
Hilbert space $\delta_0^\perp$ the operator $\cT=
P_{\delta_{0}^\perp}\cC\uphar{\delta_{0}}^\perp$. Let the operators
(matrices) $\cS:\dC\to\dC$, $\cF:\dC\to\delta^\perp_0$ and
$\cG:\delta^\perp_0\to\dC$ be given by
\[
\cS1=\bar\alpha,\quad \cF1=\begin{pmatrix}\rho_0\cr 0\cr \ldots\cr
0\end{pmatrix},\quad \cG \begin{pmatrix}h_1\cr h_2\cr \ldots\cr
h_n\cr\ldots\end{pmatrix}=\bar\alpha_1 \rho_0 h_1+\rho_1\rho_0 h_2.
\]
Hence, the matrix $\cC$ takes the block form
%\begin{equation}
%\label{block}
\[
\cC=\begin{pmatrix} \cS&\cG\cr \cF&\cT\end{pmatrix}.
%\end{equation}
\]
From \eqref{UN} it follows that
\begin{equation}
\label{FORMU}
\begin{split}
&\left\|\cG\begin{pmatrix}h_1\cr h_2\cr \ldots\cr
h_n\cr\ldots\end{pmatrix}\right\|^2=\left\|D_{\cT}\begin{pmatrix}h_1\cr
h_2\cr \ldots\cr h_n\cr\ldots\end{pmatrix}\right\|=\rho^2_0|\bar
\alpha_1 h_1+\rho_1 h_2|^2,\quad
\sD_{\cT}=\left\{\lambda(\alpha_1\delta_1+\rho_1\delta_2,\;\lambda\in\dC)\right\},\\
& \left\|\cF^*\begin{pmatrix}h_1\cr h_2\cr \ldots\cr
h_n\cr\ldots\end{pmatrix}\right\|^2=\left\|D_{\cT^*}\begin{pmatrix}h_1\cr
h_2\cr \ldots\cr
h_n\cr\ldots\end{pmatrix}\right\|=\rho^2_0|h_1|^2,\quad
\sD_{\cT^*}=\left\{\lambda\delta_1,\;\lambda\in\dC\right\},\\
&D_{\cT}
h=\rho_0(h,\alpha_1\delta_1+\rho_1\delta_2)(\alpha_1\delta_1+\rho_1\delta_2),\quad
D_{\cT^*}h=\rho_0(h,\delta_1)\delta_1,\quad h\in \ell^2(\dN)\;(\dC^N),\\
 &\cT(\alpha_1\delta_1+\rho_1\delta_2)=-\alpha_0\delta_1.
\end{split}
\end{equation}
Since $\delta_0$ is the cyclic vector for $\cC$, then by Theorem
\ref{cyclic} the unitary colligation
\begin{equation}
\label{DELTAC} \Delta_\cC=\left\{\begin{pmatrix} \cS&\cG\cr
\cF&\cT\end{pmatrix};\dC,\dC, \delta^\perp_0\right\}
\end{equation}
is prime, and $\cT$ is a completely nonunitary operator with rank
one defects on the Hilbert spaces $l_2(\dN)$ or $\dC^{N}$.

Let \be\label{C-Sh}
F(z)=\left((\cC+zI)(\cC-zI)^{-1}\delta_0,\delta_0\right),\quad
f(z)=\frac{1}{z}\frac{F(z)-1}{F(z)+1} \ee be the Carath\'eodory and
the Schur functions associated with $\cC$. By Theorems \ref{CH1} and
\ref{char} $f$ agrees with the characteristic function of $\cT$.

\begin{proposition}
\label{PROP} \begin{enumerate}
\item
 For a semi-infinite truncated CMV
matrix $\cT=\cT(\{\alpha_n\})$ the following statements are
equivalent:
\begin{enumerate}
\def\labelenumi{\rm (\roman{enumi})}
\item the matrix $\cT$ does not contain a unilateral
shift;
\item
 the matrix $\cT^*$ does not contain a unilateral
shift;
\item $\cspan\{\cT^n\delta_1,\; n=0,1,\ldots\}=\ell^2(\dN);$
\item $\cspan\{\cT^{*n}(\alpha_1\delta_1+\rho_1\delta_2),\;
n=0,1,\ldots\}=\ell^2(\dN)$;
\item
$\sum\limits_{n=0}^\infty|\alpha_n|^2=\infty;$
\item $\ln(1-|f(e^{it})|^2)\notin L^1[-\pi,\pi]$.
\end{enumerate}
\item If $\cT$ is a semi-finite truncated CMV matrix, and one of
the conditions
\begin{enumerate}
\item $\limsup_{n\to\infty}|\alpha_n|=1,$ \item
 $\lim_{n\to\infty}\alpha_n\alpha_{n+m}=0$ for
$m=1,2,\ldots$, but $\limsup_{n\to\infty}|\alpha_n|>0$
\end{enumerate}
is fulfilled, then
\[
s-\lim\limits_{n\to\infty} \cT^n=s-\lim\limits_{n\to\infty}
\cT^{*n}=0.
\]
\item
 If $\cT$ is a finite truncated CMV matrix, then
$\lim\limits_{n\to\infty}||\cT^n||=0.$
\end{enumerate}
\end{proposition}
\begin{proof}
(1) Since $\{\alpha_n\}$ are the Schur parameters of the Schur
function $f$ associated with the full CMV matrix
$\cC(\{\alpha_n\})$, and $f$ agrees with the characteristic function
of $\cT(\{\alpha_n\})$, the equivalence of the statements (a)--(f)
follows from \eqref{perp}, \eqref{shift}, \eqref{co}, \eqref{SP},
\eqref{L1}, \eqref{FORMU}, \eqref{poln}, and Theorems \ref{cyclic}
and \ref{complete}.

(2) Each condition (a) or (b) implies $f$ is inner (see subsection
4.2). Hence $\cT$ belongs to the class $C_{00}$, i.e.,
$s-\lim\limits_{n\to\infty} \cT^n=s-\lim\limits_{n\to\infty}
\cT^{*n}=0.$

(3) The function $f$ is a finite Blaschke product and so inner.
Since $\cT$ is finite-dimensional, we get
$\lim\limits_{n\to\infty}\|\cT^n\|=0.$
\end{proof}

\begin{proposition}
\label{unieq1} Let $\cT(\{\alpha_n\})$, and $\cT(\{\beta_n\})$ be
truncated CMV matrices. Then $\cT(\{\alpha_n\})$ and
$\cT(\{\beta_n\})$ are unitarily equivalent if and only if
$\beta_n=e^{it}\alpha_n$ for all $n$ and $t\in [0,2\pi)$. Moreover,
if $\cV$ is the diagonal unitary matrix of the form
\begin{equation}
\label{V} \cV={\rm{diag}}(e^{it},1,e^{it},1,\ldots),
\end{equation}
then
\begin{equation}
\label{VV} \cV\cT(\{\alpha_n\}){\cV}^{-1}=\cT(\{e^{it}\alpha_n\}).
\end{equation}
\end{proposition}
\begin{proof}  Consider two CMV matrices $\cC(\{\alpha_n\})$ and
$\cC(\{\beta_n\})$, and associated with them Schur functions
$f_\alpha$ and $f_\beta$. Since these functions agree with the
characteristic functions of $\cT(\{\alpha_n\})$ and
$\cT(\{\beta_n\})$, respectively, the operators $\cT(\{\alpha_n\})$
and $\cT(\{\beta_n\})$ are unitarily equivalent if and only if
$f_\alpha$ and $f_\beta$ differ by a scalar unimodular factor, which
in turn yields $\beta_n=e^{it}\alpha_n$ for all $n$ and $t\in
[0,2\pi)$.

Equality \eqref{VV} with $\cV$ \eqref{V} can be verified by the
direct calculation based on \eqref{cmventr1}, \eqref{cmventr2}. So
$\cT(\{\alpha_n\})$ and $\cT(\{e^{it}\alpha_n\})$ are unitarily
equivalent.
\end{proof}

\begin{remark} The similar problem for ``full'' CMV matrices
can be considered as well. Let two CMV matrices $\cC(\{\alpha_n\})$
and $\cC(\{\beta_n\})$ be unitarily equivalent by a unitary
preserving $\delta_0$. Then they are identical (see \cite[Theorem
2.3]{S3}). In general, two unitaries with simple spectra are
unitarily equivalent if and only if their spectral measures are in
the same measure class. This is a standard issue in what is called
multiplicity theory.  So, two CMV matrices are unitarily equivalent
if and only if their measures are mutually absolutely continuous.
For instance, a CMV matrix is unitarily equivalent to the free one
($\alpha_n\equiv 0$) if and only if the associated measure $\mu$ has
the property $\mu'>0$ a.e. and does not have a singular part.
\end{remark}

From \eqref{VV} it follows that
\[
\cT(\{e^{it}\alpha_n\})=e^{it\cA}\cT(\{\alpha_n\})e^{-it\cA},\;
t\in\dR,
\]
where $\cA$ is a self-adjoint diagonal matrix
$\cA=\rm{diag}(1,0,1,0\ldots)$. Hence the matrix
$\cT(\{e^{it}\alpha_n\})$ satisfies the differential equation
\[
\frac{d\cT(t)}{dt}=i\left(\cA\cT(t)-\cT(t)\cA\right)
\]
and $\cT(0)=\cT(\{\alpha_n\}).$

The next theorem states that truncated CMV matrices are models of
completely nonunitary contractions with rank one defects.
\begin{theorem}
\label{cmv1} Let $T$ be a completely nonunitary contraction with
rank one defects acting on infinite-dimensional separable Hilbert
space $H$(respectively, finite-dimensional Hilbert space). Then $T$
is unitarily equivalent to the operator acting on $l_2(\dN)$
(respectively, on $\dC^N$ in the case $\dim H=N$) determined by the
truncated CMV matrix $\cT=\cT(\{\alpha_n\})$, where $\{\alpha_n\}$
are the Schur parameters of the characteristic function of $T$. In
particular, every completely nonunitary contraction with rank one
defects is a product of co-isometric and isometric operators with
rank one defects.
\end{theorem}
\begin{proof}
Include $T$ into a prime unitary colligation
$\Delta=\left\{\begin{pmatrix} S&G\cr F&T\end{pmatrix};
\dC,\dC,H\right\}$. By Theorem \ref{cyclic} the vector $\vec
1=\begin{pmatrix}1\cr 0\end{pmatrix}$ is a cyclic for the unitary
operator $U=\begin{pmatrix} S&G\cr F&T\end{pmatrix}$. From the
results of \cite{CMV1, CMV2} (see also \cite{S1, S2}) there exists a
unique CMV matrix $\cC$ such that
\[
U=W^{-1}\cC W,\quad \delta_0=W\vec 1,
\]
where $W$ is a unitary  operator from $\dC\oplus H$ onto
$l_2(\dZ_+)$ $(\dC^{N+1})$, and $\delta_0=(1,0,0,\ldots)^t$. It
follows that the operator $W$ takes the block-operator form
\[
W=\begin{pmatrix} 1&0\cr 0&\cX \end{pmatrix}:\begin{pmatrix}\dC\cr H
\end{pmatrix}\to \begin{pmatrix}\dC\cr\delta^\perp_0
\end{pmatrix},
\]
where $\cX:H\to\delta^\perp_0$ is a unitary operator. Hence $\cT=\cX
T\cX^{-1}$, i.e., the operator $T$ is unitarily equivalent to the
operator in $\l_2(\dN)$ $(\dC^N)$ given by the truncated CMV matrix
$\cT=\cT(\{\alpha_n\})$. From representation \eqref{carat} of
$F(z)=\left((U+zI)(U-zI)^{-1}\vec 1,\vec 1\right)$ and Theorem
\ref{char} it follows that $\{\alpha_n\}$ are the Schur parameters
of the function $\ovl{\Theta_\Delta(\bar z)}$ that agrees with the
characteristic function of $T$.

Let $\cQ$ be an arbitrary unitary operator in $\delta_0^\perp$.
Since $\cT=\cL_r\cM_c$, we get
\[
T=\cX^{-1}\cT \cX=\cX^{-1}\cL_r\cM_c\cX =\cX^{-1}\cL_r
\cQ\cQ^{-1}\cM_c \cX=L\,M,
\]
where $M=\cQ^{-1}\cM_c\cX$ is an isometric operator with rank one
defect, and $L=\cX^{-1}\cL_r\cQ$ is a co-isometric operator with
rank one defect.
\end{proof}
Note that the unitary colligation \eqref{DELTAC} is unitary
equivalent to the unitary colligation \eqref{unmod}.

\subsection{The Liv$\rm{\check{s}}$ic theorem for quasi-unitary contractive
extensions and the corresponding truncated CMV matrix}

Let $V$ be an isometric operator acting on some Hilbert space $H$
with the domain $\dom V$ and the range $\ran V$. The numbers
$\dim(H\ominus\dom V)$ and $\dim (H\ominus\ran V)$ are called the
defect indices of $V$. The isometric operator $V$ is called prime if
there is no nontrivial subspace on which $V$ is unitary. In
\cite{L1, L2} M.~Liv$\rm{\check{s}}$ic developed the spectral theory
of isometric operators with equal defect indices, and their
quasi-unitary extensions. A nonunitary operator $S$ on $H$ is called
a quasi-unitary extension of the isometric operator $V$ with the
defect indices $(n,n)$, if $S$ agrees with $V$ on $\dom V$ and maps
$H\ominus \dom V$ into $H\ominus\ran V$.

Let $\vec U$ be the bilateral shift in $\ell^2(\dZ)$, i.e., $\vec
U\delta_k=\delta_{k-1}$, $k\in\dZ$, where $\{\delta_k,\; k\in\dZ\}$
is the canonical orthonormal basis in $\ell^2(\dZ)$. Define $\vec
V_0$ by
\[
\dom\vec V_0=\delta^\perp_o, \qquad \vec V_0=\vec U\uphar \dom\vec
V_0,
\]
Then $\ran \vec V_0=\delta^\perp_{-1}$. Let the quasi-unitary
extension $\vec S_0$ of $\vec V_0$ be given by $\vec S_0\delta_0=0$,
$\vec S_0\uphar\dom\vec V_0=\vec V_0$. Then each point of $\dD$ is
the eigenvalue of $\vec S_0$. So the spectrum of $\vec S_0$ agrees
with $\ovl\dD$. The following result is essentially due to M.
Liv$\rm{\check{s}}$ic \cite{L1}.
 \begin{theorem}
Let $S$ be a quasi-unitary contractive extension of a prime
isometric operator $V$ with the defect indices $(1,1)$. If the whole
open disk $\dD$ consists of the point spectrum of $S$, then $V$ and
$S$ are unitarily equivalent to $\vec V_0$ and $\vec S_{0}$,
respectively.
 \end{theorem}
Clearly, the rank of the defect operators $(I-\vec S^*_0\vec
S_0)^{1/2}$ and $(I-\vec S_0\vec S^*_0)^{1/2}$ is equal to one.
Since the point spectrum of $\vec S_0$ is $\dD$, the Sz.-Nagy--Foias
characteristic function $\Theta$ of $\vec S_0$ is identically equal
to zero. On the other hand, one can easily show (and it is well
known) that a completely nonunitary contraction with rank one
defects and zero characteristic function is unitarily equivalent to
the operator $S\oplus S^*$, where $S$ is the unilateral shift in
$\ell^2(\dN)$. So the operators $\vec S_0$ and $S\oplus S^*$ are
unitarily equivalent. Since all Schur parameters of the function
$\Theta=0$ are zeros, the corresponding truncated CMV matrix
$\cT_0=\|t_0(i,j)\|$ takes the form
\[
\cT_0=\begin{pmatrix} 0&0&0&0&0&0&\ldots\cr 0&0&0&1&0&0&\ldots\cr
1&0&0&0&0&0&\ldots\cr 0&0&0&0&0&1&\ldots\cr 0&0&1&0&0&0&\ldots \cr
\ldots&\ldots&\ldots&\ldots&\ldots&\ldots&\ldots
\end{pmatrix},
\]
i.e., $t_0(2k,2k+2)=t_0(2k+1, 2k-1)=1,$ $k\ge1$, and the rest
$t_0(i,j)=0$. The matrix $\cT_0$ is a submatrix of the free CMV
matrix $\cC_0$ corresponding to zero Schur parameters. Each point
$z$ of $\dD$ is the eigenvalue of $\cT_0$. The corresponding
eigensubspace is
\[
\sN_z=\{\lambda\, (0,1,0,z,0,z^2,0,z^3,\ldots)^t,\quad \lambda
\in\dC\}.
\]
 Hence, the spectrum
of $\cT_0$ is the closed unit disk $\ovl \dD$.

Let $\cV_0$ be the operator in $\ell^2(\dN)$
\begin{equation}
\label{VISO}
 \dom\cV_0=\ell^2(\dN)\ominus\{c\delta_2\}=\ker D_{\cT_0},\;
\cV_0=\cT_0\uphar\dom\cV_0.
\end{equation}
Then $\ran\cV_0=\ell^2(\dN)\ominus\{c\delta_1\}=\ker D_{\cT^*_0}$,
and $\cV_0$ is isometric with the defect indices $(1,1)$. The
contraction $\cT_0$ is the quasi-unitary extension of $\cV_0$ with
the zero characteristic function. Therefore, the truncated CMV
matrix $\cT_0$ is unitarily equivalent to the operator $\vec S_0$,
and by Liv$\rm{\check{s}}$ic theorem \cite{L1} the isometric
operator $\cV_0$ is unitarily equivalent to $\vec V_0$.

All other quasi-unitary contractive extensions of $\cV_0$ are given
by the truncated CMV matrices $\cT=\|t(i,j)\|$
\begin{equation}
\label{TTTO} \cT=\begin{pmatrix}0&-re^{i\f}&0&0&0&0&\ldots\cr
0&0&0&1&0&0&\ldots\cr 1&0&0&0&0&0&\ldots\cr 0&0&0&0&0&1&\ldots\cr
0&0&1&0&0&0&\ldots \cr
\ldots&\ldots&\ldots&\ldots&\ldots&\ldots&\ldots
\end{pmatrix},
\end{equation}
i.e., $t(2k,2k+2)=t(2k+1, 2k-1)=1,$ $k\ge 1,$ $t(1,2)=-re^{i\f}$,
$r\in (0,1)$, $\f$ is an arbitrary number from the interval $[0,
2\pi)$, and the rest $t(i,j)=0$. The characteristic function of
$\cT$ is the constant function $\Theta=r e^{i\f}$. The spectrum of
each such matrix is the unit circle $\dT$. Because
$|\Theta^{-1}|=r^{-1}$, each of such matrix is similar to unitary
matrix \cite [Theorem IX.1.2]{SF}.

The matrices $\cT_0$ and $\cT$ contain the shift
\[
\dom \cW=\cspan\{\delta_1,\delta_3,\ldots,\delta_{2n-1},\ldots\},\;
\cW\left(\sum_{n=1}^\infty h_n\delta_{2n-1}\right)=\sum_{n=1}^\infty
h_n\delta_{2n+1}.
\]
The matrices $\cT^*_0$ and $\cT^*$ contain the shift
\[
\dom \cW_*=\cspan\{\delta_2,\delta_4,\ldots,\delta_{2n},\ldots\},\;
\cW_*\left(\sum_{n=1}^\infty h_n\delta_{2n}\right)=\sum_{n=1}^\infty
h_n\delta_{2n+2}.
\]

Let $T$ be a completely nonunitary contraction with rank one defects
and the constant characteristic function $\Theta$,
$0<|\Theta(z)|=r<1$. Then by Theorem \ref{cmv1} $T$ is unitarily
equivalent to the truncated CMV matrices \eqref{TTTO}.

\subsection{Sub-matrices of truncated CMV matrices and iterates of their Schur
functions} Along with truncated CMV matrices $\cT(\{\alpha_n\})$
\eqref{MODEL}, we consider here truncated CMV matrices
$\wt\cT(\{\alpha_n\})$ obtained from the alternate CMV matrix
$\wt\cC(\{\alpha_n\})$ \eqref{altcmvmatr} by the same procedure. The
matrix $\wt\cT(\{\alpha_n\}$ is the transpose of $\cT(\{\alpha_n\})$
\begin{equation}
\label{TMODEL} \wt\cT=\begin{pmatrix}
-\bar{\alpha}_1\alpha_0&\bar{\alpha}_2\rho_1&\rho_2\rho_1&0&\ldots\cr
-\rho_1\alpha_0&-\bar{\alpha}_2\alpha_1&-\rho_2\alpha_1&0&\ldots\cr
0&\bar{\alpha}_3\rho_2&-\bar{\alpha}_3\alpha_2&\bar{\alpha_4}\rho_3&\ldots\cr
0&\rho_3\rho_2&-\rho_3\alpha_2&-\bar{\alpha}_4\alpha_3&\ldots\cr
\ldots&\ldots&\ldots&\ldots&\ldots
\end{pmatrix},
\end{equation}
and
\[
\wt\cT(\{\alpha_n\})=\cT^t(\{\alpha_n\})=(\cM_c)^t(\cL_r)^t=\cM_r\cL_c.
\]
As in Section 6.1, it is not hard to show that
$\wt\cT(\{\alpha_n\})$ is a completely nonunitary contraction with
rank one defects, and its characteristic function $\wt f$ agrees
with the Schur function associated with Verblunsky coefficients
(Schur parameters) $\{\alpha_n\}$. Indeed (cf. \eqref{C-Sh})
\[
(\wt\cC+zI)(\wt\cC-zI)^{-1}=(\cC^t+zI)(\cC^t-zI)^{-1}=
\left((\cC+zI)(\cC-zI)^{-1}\right)^t, \]
 and so
 $\wt F(z):=\left((\wt\cC+zI)(\wt\cC-zI)^{-1}\delta_0,\delta_0\right)=F(z)$,
$\wt f=f$, as claimed. So, the matrices $\cT(\{\alpha_n\})$ and
$\wt\cT(\{\alpha_n\})$ are unitarily equivalent.

Denote by $\cT^{(k)}$ ($\wt\cT^{(k)}$) the matrix obtained from
$\cT$ ($\wt\cT)$ by deleting the first $k$ rows and columns. The
following result provides the characteristic function of
$\cT^{(k)}$.
\begin{theorem}
\label{DER} Let $\mu$ be a probability measure on $\dT$ with
Verblunsky coefficients $\{\alpha_n\}_{n=0}^N$, $N\le\infty$, and
let $f$, $\cC(\{\alpha_n\})$, $\wt\cC(\{\alpha_n\})$,
$\cT(\{\alpha_n\})$, $\wt\cT(\{\alpha_n\})$ be the corresponding
Schur function, CMV and truncated CMV matrices, respectively. Then
$\cT^{(k)}$, $\wt\cT^{(k)}$ are completely nonunitary contractions
with rank one defects, and the following relations hold:
\[
%\begin{split}
\cT^{(2m-1)}(\{\alpha_n\}_{n=0}^N)=\wt\cT(\{\alpha_n\}_{n=2m-1}^N),
\quad
\cT^{(2m)}(\{\alpha_n\}_{n=0}^N)=\cT(\{\alpha_n\}_{n=2m}^N),\quad
m=1,2,\ldots.
%\end{split}
\]
So, the characteristic function of $\cT^{(k)}$ agrees with the
$k^{th}$ Schur iterate of $f$.
\end{theorem}
\begin{proof}
The relations
\[
\cT^{(1)}(\{\alpha_n\}_{n=0}^N)=\wt\cT(\{\alpha_n\}_{n=1}^N), \qquad
\wt\cT^{(1)}(\{\alpha_n\}_{n=1}^N)=\cT(\{\alpha_n\}_{n=2}^N)
\]
follows directly from \eqref{MODEL} and \eqref{TMODEL}. The rest is
a matter of simple induction and the definition of the $k^{th}$
Schur iterates.
\end{proof}

The relation between characteristic functions of the sub-matrices
$\cT^{(k)}(\{\alpha\}_{n=0}^N)$ and the $k^{th}$ Schur iterates
established in the above mentioned theorem is a complete analog of
the result concerning the connections between m-functions of a
Jacobi matrix and its sub-matrices \cite{GS}.

Let us now go back to the model of Section \ref{s5}.
\begin{theorem}
Let $\mu$ be a probability measure on $\dT$ with Verblunsky
coefficients $\{\alpha_n\}_{n=0}^N$, $N\le \infty$. Consider three
subspaces in $L^2(\dT,\mu)$
\[
\cH_{2m}:=\span \{1,\zeta,\bar \zeta,\zeta^2,\bar\zeta
^2,\ldots,\zeta^m,\bar\zeta^m\},
\]
\[
\cH_{2m-1}:=\span \{1,\zeta,\bar \zeta,\zeta^2,\bar\zeta
^2,\ldots,\bar\zeta^{m-1},\zeta^m\}, \quad \wt\cH_{2m-1}:=\span
\{1,\bar\zeta,
\zeta,\bar\zeta^2,\zeta^2,\ldots,\zeta^{m-1},\bar\zeta^m\}.
\]
Denote by $\sH_{2m}$ $(\sH_{2m-1},\ \wt\sH_{2m-1})$ their orthogonal
complements in $L^2(\dT,\mu)$, and by $P_{2m}$ $(P_{2m-1},\ \wt
P_{2m-1})$ the orthogonal projections onto $\sH_{2m}$ $(\sH_{2m-1},\
\wt\sH_{2m-1})$, respectively. Then the operators
\begin{equation}
\label{TK}
%\begin{split}
\sT_k h(\zeta)= P_k \left(\zeta h(\zeta)\right),\;
h(\zeta)\in\sH_k,\qquad \wt\sT_{2m-1} h(\zeta)= \wt P_m \left(\zeta
h(\zeta)\right),\; h(\zeta)\in\wt \sH_{2m-1},
%\end{split}
\end{equation}
are completely nonunitary contractions with rank one defects. The
characteristic function of  $\sT_k$ agrees with the $k^{th}$ Schur
iterate of the Schur function $f(\mu)$, the characteristic function
of $\wt\sT_{2m-1}$ with $(2m-1)^{th}$ Schur iterate of $f(\mu)$. So,
the operator $\sT_{k}$ is unitarily equivalent to the operator
\begin{equation}
\label{mmmod} h(\zeta)\to P_0^{(k)}\left(\zeta h(\zeta)\right),
\quad h(\zeta)\in
L^2\left(\dT,d\mu(\{\alpha_n\}_{n=k}^N)\right)\ominus \dC,
\end{equation}
where $P_0^{(k)}$ is the orthogonal projection onto
$L^2\left(\dT,d\mu(\{\alpha_n\}_{n=k}^N)\right)\ominus \dC$. In
addition, $\sT_{2m-1}$ is unitarily equivalent to $\wt\sT_{2m-1}$.
\end{theorem}
\begin{proof} Recall that  CMV matrices
$\cC(\{\alpha_n\}$ and $\wt\cC(\{\alpha_n\})$ represent the unitary
operator $Uh(\zeta)=\zeta h(\zeta)$ in $L^2(\dT,d\mu(\{\alpha_n\}))$
with respect to the complete orthonormal systems $\{\chi_n\}$ and
$\{x_n\}$, respectively. Moreover
\[
\begin{split}
&\cH_{2m}=
\span\{\chi_0,\chi_1,\ldots,\chi_{2m}\}=\span\{x_0,x_1,\ldots,
x_{2m}\},\\
&\cH_{2m-1}=\span\{\chi_0,\chi_1,\ldots,\chi_{2m-1}\},\\
&\wt\cH_{2m-1}=\span\{x_0,x_1,\ldots, x_{2m-1}\}.
\end{split}
\]
Since $\cT(\{\alpha_n\}_{n=0}^N\})$
($\wt\cT(\{\alpha_n\}_{n=0}^N\})$ is the matrix of $\sT$
\eqref{UNMODL2} with respect to the basis $\{\chi_n\}_{n=1}^N$
$(\{x_n\}_{n=1}^N$), the operators $\sT_{2m}$, $\sT_{2m-1}$, and
$\wt\sT_{2m-1}$ have the matrices $\cT^{(2m)}$, $\cT^{(2m-1)}$, and
$\wt\cT^{(2m-1)}$, respectively. From Theorem \ref{DER} it follows
that $\sT_k$ are completely nonunitary contractions with rank one
defects for all $k$, and their characteristic functions agree with
the $k^{th}$ Schur iterates of $f$.
 By Theorems \ref{DER} and \ref{L2} the operator $\sT_k$ is unitarily
 equivalent to the operator
 given by \eqref{mmmod}. We also have
\[
\wt\cT^{(2m-1)}(\{\alpha_n\}_{n=0}^N)=\cT(\{\alpha_n\}_{n=2m-1}^N).
\]
Therefore, the characteristic function of
$\wt\sT^{(2m-1)}(\{\alpha_n\}_{n=0}^N)$ agrees with $(2m-1)^{th}$
iterate $f_{2m-1}$ of $f$, and hence the operators
$\wt\sT^{(2m-1)}(\{\alpha_n\}_{n=0}^N)$ and
$\sT^{(2m-1)}(\{\alpha_n\}_{n=0}^N)$ are unitarily equivalent.
\end{proof}

We complete the section with the general result from the
contractions theory which is proved with the help of the truncated
CMV model.
\begin{proposition}
Let $T$ be a completely nonunitary contraction with rank one defects
in a separable Hilbert space $H$, $\dim H\ge 2$, and let $P_{\ker
D_{T^*}}$, $P_{\ker D_{T}}$ be the orthogonal projections onto $\ker
D_{T^*}$ and $\ker D_{T}$ in $H$, respectively. Then the operators
\[
 T_1:=P_{\ker D_{T^*}}\,T\uphar \ker D_{T^*},\; \wt T_1:=P_{\ker
 {D_T}}T\uphar \ker D_{T}
\]
are  unitarily equivalent completely nonunitary contractions with
rank one defects, and their characteristic functions agree with the
function
\[
h_1(z):=\frac{1}{z}\,\frac{h(z)-h(0)}{1-\overline{h(0)} h(z)},
\]
where $h$ is the characteristic function of $T$.
\end{proposition}
\begin{proof} By Theorem \ref{cmv1} the operator $T$ is unitarily
equivalent to the truncated CMV matrices
$\cT=\cT(\{\alpha_n\}_{n=0}^N)$ and
$\wt\cT=\wt\cT(\{\alpha_n\}_{n=0}^N)$, where $\{\alpha_n\}_{n=0}^N$
are the Schur parameters of $h$, $N\le\infty$. So, there exists a
unitary operators $V,\,\wt V:\delta_0^\perp\to H$ such that
 \[
V\cT V^{-1}=\wt V\wt\cT \wt V^{-1}=T.
 \]
It follows that
 \[
VD_{\cT^*}V^{-1}=D_{T^*},\qquad \wt VD_{\wt\cT}\wt V^{-1}=D_{T},
 \]
and  hence $V\ker D_{\cT^*}=\ker D_{T^*}$,  $\wt V\ker
D_{\wt\cT}=\ker D_{T}$. Due to \eqref{FORMU} we have
\[
\sD_{\cT^*}=\sD_{\wt\cT}=\span\{\delta_1\} \]
 and
\[
\cT^{(1)}=P_{\ker D_{\cT^*}}\cT\uphar\ker D_{\cT^*},\qquad \wt
\cT^{(1)}=P_{\ker D_{\wt\cT}}\wt\cT\uphar\ker D_{\wt\cT}.
\]
Hence
\[
V\cT^{(1)}V^{-1}=T_1,\qquad \wt V\wt\cT^{(1)}\wt V^{-1}=\wt T_1.
\]
Now from Theorem \ref{DER} it follows that $T_1$ and $\wt T_1$ are
completely nonunitary contractions with rank one defects, and their
characteristic functions agree with the first Schur iterate $h_1$ of
$h$. Hence $T_1$ and $\wt T_1$ are unitarily equivalent.

\end{proof}

\section{Inverse spectral problems for finite and semi-infinite truncated CMV matrices}
Consider a $N\times N$ truncated CMV matrix
\begin{equation}
\label{NMODEL} \cT=\cT (\{\alpha_n\})=\begin{pmatrix}
-\bar{\alpha}_1\alpha_0&-\rho_1\alpha_0&0&\ldots&0\cr
\bar{\alpha}_2\rho_1&-\bar{\alpha}_2\alpha_1&\bar{\alpha}_3\rho_2&\ldots&0\cr
\rho_2\rho_1&-\rho_2\alpha_1&-\bar\alpha_3\alpha_2&\ldots&0\cr
\ldots&\ldots&\ldots&\ldots&\bar\alpha_N\rho_{N-1}\cr
\ldots&\ldots&\ldots&-\rho_{N-1}\alpha_{N-2}&-\bar\alpha_N\alpha_{N-1}
\end{pmatrix}
\end{equation}
(for even $N$ it looks a bit different). The problem under
investigation in the present section is the reconstruction of the
matrix $\cT$ \eqref{NMODEL} from either the complete set of its
eigenvalues or from the mixed spectral data: the part of the
spectrum and the part of the parameters $\alpha_n(\cT)$.

\subsection{Existence of a finite truncated CMV matrix
 with the given spectrum}
\begin{theorem}
\label{InvSp} Let $z_1, z_2,\ldots, z_N$ be not necessarily distinct
numbers from the open unit disk. Then there exists a truncated
$N\times N$ CMV matrix $\cT$ \eqref{NMODEL} which has eigenvalues
$z_1, z_2,\dots,z_N$, counting their algebraic multiplicities. Such
matrix is determined uniquely up to multiplication of its parameters
$\alpha_n(\cT)$ by the same unimodular factor.
\end{theorem}
\begin{proof}
Let
 \be\label{blash} b(z)=e^{i\f}\prod\limits_{k=1}^N\frac{z-{z}_k}{1-\bar z_k\,
z},\quad z\in\dD,\quad \f\in [0,2\pi). \ee We want to show that $b$
is the characteristic function of a truncated CMV matrix $\cT$
(\ref{NMODEL}). Put
\[
F(z)=\frac{1+zb(z)}{1-zb(z)}\,,
\]
which is a rational function with $N+1$ distinct simple poles lying
on $\dT$, $\RE F(z)>0$, $z\in\dD$, and $ F(0)=1$. It follows that
there exists a probability measure $d\mu$ on the unit circle
supported at those poles, so that
\[
F(z)=\int_{\dT}\frac{\zeta+z}{\zeta-z}\,d\mu(\zeta).
\]
Let $\{\alpha_0,\ldots,\alpha_{N-1},\alpha_N\}$ be the Schur
parameters of $b$, that is the same as the Verblunsky coefficients
of $\mu$. Construct the $(N+1)\times (N+1)$ unitary CMV matrix $\cC
$ of the form \eqref{cmvmatr}. Then
\[
F(z)=\left((\cC+zI)(\cC-zI)^{-1}\delta_0,\delta_0\right),\;|z|<1,
\]
where $\delta_0=(1,0,\ldots,0)^t\in\dC^{N+1}$. Let $\cT$ be $N\times
N$ be truncated CMV matrix of the form \eqref{NMODEL}. $\cC$ has the
block form
\[
\cC=\begin{pmatrix}\cS& \cG\cr  \cF&\cT\end{pmatrix},
\]
where $ \cS=\bar\alpha_0$, $
\cG=\begin{pmatrix}\bar\alpha_1\rho_0,&\rho_1\rho_0,&0,&\ldots,&0\end{pmatrix}$,
and $\cF=\begin{pmatrix}\rho_0\cr 0\cr \ldots\cr 0\end{pmatrix}$.
Since $\delta_0$ is a cyclic vector for $\cC$, the unitary
colligation $\Delta=\left\{\begin{pmatrix} \cS& \cG\cr
\cF&\cT\end{pmatrix},\;\dC,\dC,\dC^{N}\right\}$ is prime. Hence
$\cT$ is a completely nonunitary contraction with rank one defect
operators. Let $\Theta_{\Delta}(z)$ be the transfer function of
$\Delta$. By Theorem \ref{char} we have
\[
\ovl{\Theta_{\Delta}(\bar z)}=\frac{1}{z}\frac{F(z)-1}{F(z)+1}\,,
\qquad \Theta_\Delta(z)=\ovl{b(\bar z)}.
\]
So $b(z)$ agrees with the characteristic function of $\cT$.
Therefore $\cT$ has eigenvalues $z_1,\ldots, z_N$, counting their
algebraic multiplicities \cite{SF}.

Finally, let $\cT(\{\alpha_n\})$ and $\cT(\{\beta_n\})$ be two such
matrices. Each of them is a completely nonunitary matrix with rank
one defects, and their characteristic functions agree with $b$
\eqref{blash}. Hence they are unitarily equivalent, and Proposition
\ref{unieq1} completes the proof.
\end{proof}
\begin{example} Let $T$ be a completely nonunitary contraction with rank
one defects on $N$-dimensional Hilbert space, and let $T$ have just
one eigenvalue $z=0$ of the algebraic multiplicity $N$. Then its
characteristic function agrees with $f(z)=e^{i\f}z^N$. The
corresponding Schur parameters are
$\{\underbrace{0,\dots,0}_{N},e^{i\f}\}.$ It follows that $\rho_n=1$
for $n=0,\ldots, N-1$. Hence $T$ is unitarily equivalent to the
$N\times $N truncated CMV matrix $\cT_N$ (see the expressions for
$\cT_5$ and $\cT_6$):
\[
\cT_5=\begin{pmatrix} 0&0&0&0&0\cr 0&0&0&1&0\cr 1&0&0&0&0\cr
0&0&0&0&e^{i\f}\cr 0&0&1&0&0\end{pmatrix},\; \cT_6=\begin{pmatrix}
0&0&0&0&0&0\cr 0&0&0&1&0&0\cr 1&0&0&0&0&0\cr 0&0&0&0&0&1\cr
0&0&1&0&0&0\cr 0&0&0&0&e^{i\f}&0\end{pmatrix}.
\]
\end{example}
\subsection{Uniqueness and reconstruction of a finite truncated CMV matrix from
mixed spectral data}

It is easily seen from \eqref{NMODEL} that a truncated $N\times N$
CMV matrix $\cT$ is completely determined by $N+1$ independent
parameters $\alpha_j(\cT)$, $j=0,1,\ldots,N$. The problem we discuss
here is whether $\cT$ can be restored from the part of its spectrum
(the eigenvalues $z_1,\ldots,z_m$, of the algebraic multiplicity
$l_k$, $k=1,\ldots,m$, with $l_1+\ldots+l_m=r$), and the {\it first}
$N-r+1$ parameters $\alpha_0(\cT),\ldots,\alpha_{N-r}(\cT)$. As we
will see later on, the solution of this problem is unique (if it
exists).

We begin with a simple result from complex analysis. We don't know
where exactly it appears in the literature, but by all means it is
known to experts.

\begin{lemma}\label{np}
Let $z_1,\ldots,z_m$ be distinct points in $\dD$, $l_1,\ldots,l_m$
positive integers, and $r=l_1+\ldots+l_m$. Suppose that the
Nevanlinna-Pick interpolation problem with multiple nodes
\be\label{np1} b^{(j)}(z_k)=w_k^{(j)}, \qquad j=0,1,\ldots,l_k-1,
\quad k=1,2,\ldots,m \ee has two solutions $b_1$ and $b_2$, both the
Blaschke products of order $\leq r-1$. Then $b_1=b_2$.
\end{lemma}
\begin{proof}
Assume first that $z_k\not=0$, $w_k^{(0)}\not=0$, $k=1,\ldots,m$.
Given a Blaschke product $s$, we see by differentiating the equality
$\ovl{s(1/\bar z)}=s^{-1}(z)$ that
\[
\ovl{s^{(j)}\left(\frac1{\bar
z}\right)}=\frac{P_j(s(z),s'(z),\ldots,s^{(j)}(z))}{s^{2^j}(z)}\,,
\]
where $P_j$ is a polynomial of its variables. Hence
\[
\ovl{s^{(j)}\left(\frac1{\bar
z_k}\right)}=\frac{P_j(s(z_k),\ldots,s^{(j)}(z_k))}{s^{2^j}(z_k)}\,,
\qquad k=1,2,\ldots,m
\]
so we have
\[
b_1^{(j)}(z_k)=b_2^{(j)}(z_k), \quad b_1^{(j)}\left(\frac1{\bar
z_k}\right)=b_2^{(j)}\left(\frac1{\bar z_k}\right), \qquad
j=0,1,\ldots,l_k-1, \quad k=1,2,\ldots,m.
\]
Then for the difference $u=b_1-b_2$ the relations \be\label{np2}
u^{(j)}(z_k)=u^{(j)}\left(\frac1{\bar z_k}\right)=0, \qquad
j=0,1,\ldots,l_k-1, \quad k=1,2,\ldots,m. \ee hold. Let now
\[
b_l(z)=\frac{p_l(z)}{q_l(z)}\,, \quad l=1,2, \qquad
u(z)=\frac{p_1(z)q_2(z)-p_2(z)q_1(z)}{q_1(z)q_2(z)}=\frac{p(z)}{q(z)}\,,
\]
where $p,q$ are polynomials of degree $\leq 2r-2$. The Leibniz
formula
\[
u^{(n)}(z)=\sum_{k=0}^n
\frac{n!}{k!(n-k)!}\,p^{(k)}(z)\left(\frac1{q}\right)^{(n-k)}(z)
\]
shows by induction that \eqref{np2} imply \be\label{np3}
p^{(j)}(z_k)=p^{(j)}\left(\frac1{\bar z_k}\right)=0, \qquad
j=0,1,\ldots,l_k-1, \quad k=1,2,\ldots,m. \ee But $\deg p\leq 2r-2$,
and there are $2r$ conditions in \eqref{np3}, so $p\equiv 0$, as
needed.

Assume next that $z_k\not=0$, $k=1,\ldots,m$ and some of $w_k^{(0)}$
are zero. Take $\varepsilon\in\dD$, $\varepsilon\not=w_k^{(0)}$ and
put
\[
s_0:=\frac{z-\varepsilon}{1-\bar\varepsilon z}\,, \qquad \wh
b_l(z):=s_0(b_l(z)), \qquad l=1,2.
\]
Then both $\wh b_1$ and $\wh b_2$ are Blaschke products of order
$\leq r-1$ which solve the interpolation problem
\[
\wh b_l^{(j)}(z_k)=\wh w_k^{(j)}, \qquad j=0,1,\ldots,l_k-1, \quad
k=1,2,\ldots,m, \quad l=1,2,
\]
where $\wh w_k^{(0)}=s_0(w_k^{(0)})\not=0$ and $\wh w_k^{(j)}=\left(
s_0(b_l(z))\right)^{(j)}\Bigr|_{z=z_k}$. The above argument applied
to $\wh b_l$ gives $\wh b_1=\wh b_2 \Rightarrow b_1=b_2$, as needed.

Finally, assume that $z_1=0$. Let $\varepsilon\not=-z_k$ for all
$k$, and put
\[
\wt b_l(z):=b_l(s_0(z)), \qquad l=1,2.
\]
Then the Blaschke products $\wt b_1$, $\wt b_2$ of order $\leq r-1$
satisfy
\[
\wt b_l^{(j)}(\wt z_k)=\wt w_k^{(j)}, \qquad j=0,1,\ldots,l_k-1,
\quad k=1,2,\ldots,m, \quad l=1,2
\]
and $\wt z_k=(z_k+\varepsilon)(1+\bar\varepsilon z_k)^{-1}\not=0$.
Hence $\wt b_1=\wt b_2$, and so $b_1=b_2$. The proof is complete.

\end{proof}

\begin{theorem}
\label{mix1} Let $z_1,\ldots,z_m$ be distinct nonzero points in
$\dD$, $l_1,\ldots,l_m$ be positive integers, and
$r=l_1+\ldots+l_m\leq N$. Let $\alpha_0,\ldots,\alpha_{N-r}\in\dD$.
If there exists a $N\times N$ truncated CMV matrix $\cT$
\eqref{NMODEL} such that $z_1,\ldots z_m$ are eigenvalues of $\cT$
with the algebraic multiplicities $l_1,\ldots, l_m$, and
$\alpha_j(\cT)=\alpha_j$, $j=0,\ldots,N-r$, then this matrix is
unique.
\end{theorem}
\begin{proof}
If the required $\cT$ exists then its characteristic function
$\Theta_{\cT}(z)$ is the Blaschke product of order $N$ and of the
form
\begin{equation}
\label{one} b(z)=e^{it}\prod\limits_{k=1}^m\left(\frac{z-z_k}{1-\bar
z_k z}\right)^{l_k}\,\prod\limits_{j=1}^{N-r}\frac{z-v_j}{1-\bar v_j
z},
\end{equation}
with the given first $N-r+1$ Schur parameters
$\alpha_0(b),\ldots,\alpha_{N-r}(b)$. Our goal is to prove the
uniqueness of such function $b$.

According to the result of Schur \cite{Schur} (see Section 4.2) the
set of all Schur functions $b$ with given first $N-r+1$ Schur
parameters is parametrized by \be\label{shu}
b(z)=\frac{A(z)+zB^*(z)s(z)}{B(z)+zA^*(z)s(z)}\,, \ee where $s(z)$
is an arbitrary Schur function, and $A,B$ are polynomials of degree
at most $N-r$. Since $b$ is the Blaschke product of order $N$, it is
clear that so is $s(z)$, $\deg s(z)=r-1$, and
\[
\cS b=\{\alpha_0,\ldots,\alpha_{N-r},\alpha_0(s),\ldots,
\alpha_{r-1}(s)\}.
\]
Let us solve \eqref{shu} for $s$:
\[
s(z)=\frac{A(z)-B(z)b(z)}{-zB^*(z)+zA^*(z)b(z)}\,,
\]
so $s(z)$ satisfies the Nevanlinna-Pick interpolation problem
\eqref{np1}, where $w^{(j)}_k$ are completely determined from the
given nonzero $z_k$'s and $\alpha_j$'s. By Lemma \ref{np} there is
at most one such $s(z)$, and the uniqueness of $b$ is proved.
\end{proof}

\begin{remark}
Suppose that $z_1,\ldots,z_m$ are distinct nonzero points in $\dD$,
and $l_1+\ldots+l_m=N$, so the only $\alpha_0$ is prescribed. It is
clear that $\alpha_0$ is completely determined by the choice of
$z_j$ and their multiplicities $l_j$:
\[
b(z)=e^{it}\prod\limits_{k=1}^m\left(\frac{z-z_k}{1-\bar z_k
z}\right)^{l_k}, \quad
\alpha_0=b(0)=e^{it}\prod_{j=1}^m(-z_k^{l_k}).
\]
So for all other $\alpha_0$ the inverse problem has no solution.
\end{remark}

In the case when one of the eigenvalues is zero, all three
possibilities (no solution, unique solution, and infinitely many
solutions) may occur for the inverse problem in question. For
instance, there is no solution at all as long as $z_1=0$,
$\alpha_0\not=0$. Assume next, that $r=l_1=1$, $z_1=0$, and the
points $\alpha_0,\alpha_1,\ldots,\alpha_{N-1}$ are taken in $\dD$,
with the only restriction $\alpha_0=0$, $\alpha_1\not=0$. The
Blaschke products $b_\gamma$ with the Schur parameters
$\{\alpha_0,\alpha_1,\ldots,\alpha_{N-1};\gamma\}$ and arbitrary
$\gamma\in\dT$ are of the form
\[
b_\gamma(z)=e^{it}z\prod_{j=1}^{N-1}\frac{z-v_j}{1-\bar v_j z},
\]
and the corresponding $N\times N$ truncated CMV matrices
$\cT_\gamma$ solve the problem.

Finally, assume that except for the zero eigenvalue of multiplicity
$k$ ($z_1=z_2=\ldots=z_k=0$), a few more nonzero (and not
necessarily distinct) eigenvalues $\lambda_1,\ldots,\lambda_r$ are
given, as well as the points $\alpha_0=\ldots=\alpha_{k-1}=0$,
$\alpha_k\not=0,\ldots,\alpha_{N-r}$ in $\dD$. If the solution of
the corresponding mixed inverse problem $\cT$ exists, its
characteristic function takes the form
\[
b(z)=e^{it}z^k\prod\limits_{j=1}^r\frac{z-\lambda_j}{1-\bar
\lambda_j z}\,g(z),
\]
where $g$ is the Blaschke product of order $N-k-r$, $g(0)\not=0$,
and the first $N-k-r+1$ Schur parameters of $h=z^{-k}b$ are given
numbers $\alpha_k,\ldots,\alpha_{N-r}$\ . Clearly, $h$ is exactly
the $k^{th}$ Schur iterate of $b$. If the required truncated CMV
matrix $\cT$ exists, then by Theorem \ref{DER} the characteristic
function of $\cT^{(k)}$ agrees with $h$. It follows now from Theorem
\ref{mix1} that $\cT^{(k)}$ is unique, and since $\alpha_j(\cT)=0$,
$j=0,\ldots,k-1$, the matrix $\cT$ is unique as well.

The situation changes dramatically if we assume that the {\it last}
parameters of $\cT$ \eqref{NMODEL} are known. In this case we can
prove the existence, but not the uniqueness of the solution.

\begin{theorem}
\label{mix2} Let $z_1, \ldots, z_m$ and
$\alpha_m,\ldots\alpha_{N-1}$ be two collections of arbitrary
complex number from the open unit disk, and let $\alpha_N\in\dT$.
Then there exists a $N\times N$ truncated CMV matrix $\cT$ of the
form \eqref{NMODEL} such that
\begin{enumerate}
\def\labelenumi{\rm (\roman{enumi})}
\item $z_1,\ldots,z_m$ are eigenvalues of $\cT$, counting the
algebraic multiplicity, \item $\alpha_n(\cT)=\alpha_n$,
$n=m,m+1,\ldots,N$.
\end{enumerate}
\end{theorem}
\begin{proof} By Theorem \ref{TTT} there exists a Blaschke product
$b(z)$ of order $N$ such that $b(z_k)=0$, $k=1,\ldots,m$, with the
Schur parameters
\[
\alpha_n(b)=\alpha_n, \qquad n=m,m+1,\ldots,N.
\]
Take now the matrix $\cT$ \eqref{NMODEL} with
$\alpha_n(\cT)=\alpha_n$, $n=0,1,\ldots,N$. By Theorem \ref{char}
the characteristic function of $\cT$ agrees with $b(z)$, that
completes the proof.
\end{proof}
Theorem 7.6 thereby says that a $N\times N$ truncated CMV matrix
$\cT$ can be reconstructed from its $m$ eigenvalues and the lower
principal block of order $N-m$. The latter is either the truncated
CMV matrix $\cT(\{\alpha_n\}_{n=m}^N)$ or its transpose $\wt\cT$.

\subsection{Inverse problem for semi-infinite truncated CMV matrix}
 In this subsection we consider the criterion when given
  complex numbers $z_n$, $n=1,2,\ldots$ from $\dD$ are the
  eigenvalues counting algebraic multiplicity of some semi-infinite
  truncated CMV matrix.
  \begin{proposition}
   Given complex numbers $z_n$, $n=1,2,...$ are eigenvalues counting
   algebraic multiplicity of some semi-infinite truncated CMV matrix
   if and only if
 \[
 \sum\limits_{n=1}^{\infty}(1-|z_n|)< \infty.
 \]
 \end{proposition}
\begin{proof}
 The convergence of the sum is equivalent to the convergence of the Blashke
 product
$$ b(z)=\prod\limits_{k=1}^{\infty}\frac{\bar z_k}{z_k}\,\frac{z_k-z}{1-\bar z_k
z}.
$$
 Let $\{\alpha_n\}$ be the Schur parameters of $b$. The
characteristic function of the truncated CMV matrix
$\cT(\{\alpha_n\})$ agrees with $b$. Hence the eigenvalues of
$\cT(\{\alpha_n\})$ are precisely the complex numbers $\{z_n\}$.
\end{proof}

\end{document}